\documentclass[runningheads]{llncs}
\usepackage{isabelle_indepCH,isabellesym_indepCH}
\usepackage{booktabs,array,threeparttable}
\usepackage{amsmath,amsfonts,amssymb}
\usepackage{bbm}  %
\usepackage{tikz}
\usepackage[normalem]{ulem}
\usepackage[english]{babel}
\usepackage{multidef}
\usepackage{verbatim}
\usepackage{stmaryrd} %
\usepackage{hyperref}
\usepackage{xcolor}
\usepackage{framed}
\usepackage[numbers]{natbib}
\usepackage{relsize}
\DeclareRobustCommand{\isactrlbsub}{\emph\bgroup\math{}\sb\bgroup\mbox\bgroup\isaspacing\itshape\smaller}
\DeclareRobustCommand{\isactrlesub}{\egroup\egroup\endmath\egroup}
\DeclareRobustCommand{\isactrlbsup}{\emph\bgroup\math{}\sp\bgroup\mbox\bgroup\isaspacing\itshape\smaller}
\DeclareRobustCommand{\isactrlesup}{\egroup\egroup\endmath\egroup}

\newcommand{\limp}{\longrightarrow}
\newcommand{\lsii}{\longleftrightarrow}

\DeclareMathOperator{\cf}{cf}

\DeclareMathOperator{\Fn}{Fn}

\DeclareMathOperator{\Con}{Con}

\newcommand{\modelo}[1]{\mathbf{#1}}
\newcommand{\axiomas}[1]{\mathit{#1}}
\newcommand{\clase}[1]{\mathsf{#1}}

\multidef{\clase{#1}}{Card,HC,HF,Lim,On->Ord,Reg,WF,Ord}

\DeclareMathAlphabet{\mathbbm}{U}{bbm}{m}{n} 
\newcommand{\1}{\isasymone}
\newcommand{\PP}{\mathbbm{P}}

\multidef[prefix=cal]{\mathcal{#1}}{A-Z}
\multidef{\modelo{#1}}{A,BB->B,CC->C,NN->N,QQ->Q,RR->R,ZZ->Z}

\multidef[prefix=p]{\mathbb{#1}}{A-Z}
\renewcommand{\emptyset}{\varnothing}

\newcommand{\Pow}{\mathop{\mathcal{P}}}

\newcommand{\lb}{\langle}
\newcommand{\rb}{\rangle}

\renewcommand{\phi}{\varphi}

\newcommand{\del}{\delta}

\newcommand{\ga}{\gamma}

\newcommand{\defi}{\mathrel{\mathop:}=}
\newcommand{\forces}{\Vdash}

\newcommand{\om}{\ensuremath{\omega}}

\newcommand*{\ale}[1]{\aleph_{#1}}

\newcommand{\ZC}{\axiomas{ZC}}
\newcommand{\AC}{\axiomas{AC}}
\newcommand{\DC}{\axiomas{DC}}

\newcommand{\CH}{\axiomas{CH}}
\newcommand{\ZFC}{\axiomas{ZFC}}
\newcommand{\ZF}{\axiomas{ZF}}

\newcommand{\union}{\mathop{\textstyle\bigcup}}
\newcommand{\sm}{\setminus}
\newcommand{\sbq}{\subseteq}

\newcommand{\dimg}{\text{\textup{``}}} %
\newcommand{\quine}[1]{\isa{{\isasymcdot}}#1\isa{{\isasymcdot}}}

\DeclareMathOperator{\val}{\mathit{val}}

\newcommand{\formula}{\isa{formula}}
\newcommand{\forceisa}{\mathop{\mathtt{forces}}}

\newcommand{\uscore}{\isacharunderscore}

\newcommand{\frcat}{\mathop{\isa{frc{\uscore}at}}}

\renewcommand{\isasymtturnstile}{\isamath{\Vdash}}
\renewcommand{\isacharminus}{-}
\newcommand{\session}[1]{\textit{#1}}
\newcommand{\theory}[1]{\texttt{#1}}
\renewcommand{\locale}[1]{\isa{#1}}
\newcommand{\punto}[1]{\overset{\tikz\draw[fill=black] (0,0) circle (0.6pt);}{#1}}
\newcommand{\Meta}{\mathcal{M}}
\newcommand{\type}[1]{\isa{#1}}
\newcommand{\tyi}{\type{i}}
\newcommand{\tyo}{\type{o}}
\newcommand{\prop}{\type{prop}}
\newcommand{\fun}{\Rightarrow}
\newcommand{\ALL}{\textstyle\bigwedge}

\makeatletter
\@ifpackageloaded{txfonts}\@tempswafalse\@tempswatrue
\if@tempswa
  \DeclareFontFamily{U}{txsymbols}{}
  \DeclareFontFamily{U}{txAMSb}{}
  \DeclareSymbolFont{txsymbols}{OMS}{txsy}{m}{n}
  \SetSymbolFont{txsymbols}{bold}{OMS}{txsy}{bx}{n}
  \DeclareFontSubstitution{OMS}{txsy}{m}{n}
  \DeclareSymbolFont{txAMSb}{U}{txsyb}{m}{n}
  \SetSymbolFont{txAMSb}{bold}{U}{txsyb}{bx}{n}
  \DeclareFontSubstitution{U}{txsyb}{m}{n}
  \DeclareMathSymbol{\aleph}{\mathord}{txsymbols}{64}
  \DeclareMathSymbol{\beth}{\mathord}{txAMSb}{105}
  \DeclareMathSymbol{\gimel}{\mathord}{txAMSb}{106}
  \DeclareMathSymbol{\daleth}{\mathord}{txAMSb}{107}
\fi
\makeatother

\newcommand{\quantRel}[3]{#1 #2\kern -1pt[#3]}

\newif\ifarXiv
\newif\ifIEEE

\usepackage{graphicx}
\hypersetup{
  colorlinks,
  urlcolor={blue},
  linkcolor={blue!50!black},
  citecolor={blue!50!black},
}
\DeclareUnicodeCharacter{2200}{\ensuremath{\forall}}
\DeclareUnicodeCharacter{2203}{\ensuremath{\exists}}
\DeclareUnicodeCharacter{2227}{\ensuremath{\wedge}}
\DeclareUnicodeCharacter{2228}{\ensuremath{\vee}}
\DeclareUnicodeCharacter{22C0}{\ensuremath{\bigwedge}}
\DeclareUnicodeCharacter{2261}{\ensuremath{\equiv}}
\DeclareUnicodeCharacter{27F9}{\ensuremath{\Longrightarrow}}
\DeclareUnicodeCharacter{27F7}{\ensuremath{\leftrightarrow}}
\DeclareUnicodeCharacter{27F6}{\ensuremath{\longrightarrow}}
\DeclareUnicodeCharacter{03BB}{\ensuremath{\lambda}}

\DeclareUnicodeCharacter{2208}{\ensuremath{\in}}
\DeclareUnicodeCharacter{2209}{\ensuremath{\not\in}}
\DeclareUnicodeCharacter{2286}{\ensuremath{\subseteq}}
\DeclareUnicodeCharacter{22C3}{\ensuremath{\bigcup}}
\DeclareUnicodeCharacter{21D2}{\ensuremath{\Rightarrow}}
\usepackage{array}
\newcounter{replInstCount}
\newcounter{LamReplCount}
\begin{document}
\title{The formal verification of the ctm approach to forcing%
  \thanks{Supported by Secyt-UNC projects 33620180100855CB and 33620180100465CB, and Conicet.}%
}
\titlerunning{Formalization of ctm forcing}%
\author{Emmanuel Gunther\inst{1} \and
Miguel Pagano\inst{1} \and \\
Pedro Sánchez Terraf\inst{1,2}%
\and
Matías Steinberg\inst{1}
}
\authorrunning{E.~Gunther, M.~Pagano, P.~Sánchez Terraf, M.~Steinberg}
\institute{Universidad Nacional de C\'ordoba. 
  \\  Facultad de Matem\'atica, Astronom\'{\i}a,  F\'{\i}sica y
  Computaci\'on.
  \and
    Centro de Investigaci\'on y Estudios de Matem\'atica (CIEM-FaMAF),
    Conicet. C\'ordoba. Argentina. \\
    \email{\{gunther,sterraf\}@famaf.unc.edu.ar\\
        \{miguel.pagano,matias.steinberg\}@unc.edu.ar}
}
\maketitle              %
\begin{abstract}
  We discuss some highlights of our computer-verified
  proof of the construction, given a countable transitive set-model $M$
  of $\ZFC$, of generic extensions  satisfying $\ZFC+\neg\CH$ and $\ZFC+\CH$.
  Moreover, let $\calR$ be the set of instances of the Axiom of
  Replacement. We isolated a 21-element subset $\Omega\sbq\calR$ and
  defined $\calF:\calR\to\calR$
  such that for every $\Phi\sbq\calR$ and $M$-generic $G$,
  $M\models \ZC \cup \calF\dimg\Phi \cup \Omega$ implies
  $M[G]\models \ZC \cup \Phi \cup \{ \neg \CH \}$, where $\ZC$ is
  Zermelo set theory with Choice.

  To achieve this, we worked in the proof assistant \emph{Isabelle},
  basing our development on the Isabelle/ZF library by L.~Paulson and
  others.

\end{abstract}
\section{Introduction}
\label{sec:introduction}

This paper is the culmination of our project on the computerized
formalization of the undecidability of the Continuum Hypothesis
($\CH$) from Zermelo-Fraenkel set theory with Choice ($\ZFC$), under the
assumption of the existence of a countable transitive model (ctm) of
$\ZFC$. In contrast to our reports of the previous steps towards this
goal
\cite{2018arXiv180705174G,2019arXiv190103313G,2020arXiv200109715G}, we
intend here to present our development to the mathematical logic
community. For this reason, we start with a general discussion around
the formalization of mathematics.

\subsection{Formalized mathematics}
The use of computers to assist the creation and verification of
mathematics has seen a steady grow. But the general awareness on the
matter still seems to be a bit scant (even among mathematicians
involved in foundations), and the venues devoted to the communication
of formalized mathematics are, mainly, computer science journals and
conferences: JAR, ITP, IJCAR, CPP, CICM, and others.

Nevertheless, the discussion about the subject in central mathematical
circles is increasing; there were some hints on the ICM2018 panel on
“machine-assisted” proofs
\cite{https://doi.org/10.48550/arxiv.1809.08062} and a lively
promotion by Kevin Buzzard, during his ICM2022 special plenary lecture
\cite{2021arXiv211211598B}.

Before we start an in-depth discussion, a point should be made clear:
A formalized proof is not the same as an \emph{automatic proof}. The
reader surely understands that, aside from results of a very specific sort, no current
technology allows us to write a reasonably complex (and correct)
theorem statement in a computer and expect to obtain a proof after hitting “Enter”, at
least not after a humanly feasible wait. On the other hand, it is
quite possible that the same reader has some mental image that
formalizing a proof requires making each application of Modus Ponens
explicit.

The fact is that \emph{proof assistants} are designed for the human prover to
be able to decompose a statement to be proved into smaller subgoals
which can actually be fed into some automatic tool. The balance between
what these tools are able to handle is not  easily appreciated by
intuition: Sometimes, ``trivial'' steps are not solved by them, which
can result in obvious frustration; but they would quickly solve some
goals that do not look like a ``mere computation.''

To appreciate the extent of mathematics formalizable, it is convenient to recall
some major successful projects, such as the Four Color Theorem
\cite{MR2463991}, the Odd Order Theorem
\cite{10.1007/978-3-642-39634-2_14}, and the proof the Kepler's
Conjecture \cite{MR3659768}. There is a vast mathematical corpus at
the Archive of Formal Proofs (AFP) based on Isabelle; and formalizations of
brand new mathematics like the Liquid Tensor Experiment
\cite{LTE2020,LTE2021,nature-lte} and the definition of perfectoid spaces \cite{10.1145/3372885.3373830}
have been achieved using Lean.

We will continue our description of proof assistants in
Section~\ref{sec:proof-assist-isabelle}. We kindly invite the reader
to enrich the previous exposition by reading the apt summary by
A.~Koutsoukou-Argyraki \cite{angeliki} and the interviews
therein; some of the experts consulted have also discussed
in \cite{2022arXiv220704779B} the status of formalized versus standard
proof in mathematics.

\subsection{Our achievements}
We formalized a model-theoretic rendition of forcing (Sect.~\ref{sec:forcing}), showing how to
construct proper extensions of ctms of $\ZF$ (respectively, with
$\AC$), and we formalized the basic forcing notions required to obtain
ctms of $\ZFC + \neg\CH$ and of $\ZFC + \CH$ (Sect.~\ref{sec:models-ch-negation}). No metatheoretic issues
(consistency, FOL calculi, etc) were formalized, so we were mainly
concerned with the mathematics of forcing. Nevertheless, by inspecting
the foundations underlying our proof assistant Isabelle
(Section~\ref{sec:isabelle-metalogic-meta}) it can be stated that our
formalization is a bona fide proof in $\ZF$ of the previous
constructions.

In order to reach our goals, we provided basic results that were
missing from Isabelle's $\ZF$ library, starting from ones
involving cardinal successors, countable sets, etc.
(Section~\ref{sec:extension-isabellezf}). We also extended the treatment of relativization of
set-theoretical concepts (Section~\ref{sec:tools-relativization}).

One added value that is obtained from the present formalization is
that we identified a handful of instances of Replacement which are
sufficient to set the forcing machinery up
(Section~\ref{sec:repl-instances}), on the basis of Zermelo set theory.
The eagerness to obtain this level of detail might be a consequence of
“an unnatural tendency to investigate, for the most part, trivial
minutiae of the formalism” on our part, as it was put by Cohen
\cite{zbMATH02012060}, but we would rather say that we were driven by
curiosity.

The code of our formalization can be accessed at the
AFP site, via the following link:
\begin{center}
  \url{https://www.isa-afp.org/entries/Independence_CH.html}
\end{center}

\section{Proof assistants and Isabelle/ZF}
\label{sec:proof-assist-isabelle}

Let us briefly introduce Isabelle \cite{DBLP:books/sp/Paulson94} in the large landscape of proof
assistants (“assistants” for short; also known as “interactive theorem provers”); we refer to the
excellent chapter by \citet{DBLP:series/hhl/HarrisonUW14} for a more
thorough reconstruction of the history of assistants.

It is expected that an assistant aids the human user while mechanizing
some piece of mathematics; the interaction varies from system to
system, but a common interface consists of a display showing the
current goal and assumptions. The user instructs the assistant to
modify them by means of \emph{tactics}; a %
proof is completed when the (current) goal is an instance of one of
the assumptions.

In that dialog, the user produces a script of tactics that can be
later reproduced step-by-step by the system (to check, for example,
that an imported theory is correct) or by the user to understand
the proof.

To have any value at all, the system should only allow the application
of sound tactics.
Edinburgh LCF~\cite{Gordon1979-qm} was an
influential proof assistant in which the critical code (that
constructs proofs in response to user scripts or other modules) was reduced to a small
\emph{kernel}. Hence, by verifying the correctness of the kernel, one
achieves confidence on the whole system.

The metalogic of Isabelle, as well as that of LCF, is based on higher-order logic.
In contrast, some of the other prominent assistants of today are
based on some (dependent) type theory. Both Coq~\cite{coq} and
Lean~\cite{DBLP:conf/cade/Moura021} are based on the Calculus of
Inductive Constructions
\cite{DBLP:conf/mfps/PfenningP89,DBLP:journals/iandc/CoquandH88},
while Agda~\cite{agda} is an extension of Martin-Löf type
theory~\cite{DBLP:books/daglib/0000395}. Mizar~\cite{mizar} is the
oldest assistant still used today but is far away both in terms of
foundations and architecture from Isabelle; Mizar inspired, though, the
Isar \cite{DBLP:conf/tphol/Wenzel99} dialect used in Isabelle
nowadays, which aims at the production of proof scripts that are
closer to mathematical texts.%
\footnote{%
We recommend the survey \cite{2020arXiv200909541A} by
J.~Avigad for details about the different logical foundations on which
assistants are based.}

Isabelle also inherited from LCF the possibility for the user to define
tactics to encapsulate common patterns used to solve goals. In
fact, this is the origin of the ML family of languages: a
Meta-Language for programming tactics.
In the case of Isabelle, \emph{Standard} ML is the first of the four layers
on which we worked in this assistant. Both the kernel and the automation
of proofs are coded in ML, sometimes as a substitute for induction on
formulas, as the next section explains.

\subsection{Isabelle metalogic $\Meta$}
\label{sec:isabelle-metalogic-meta}

The second layer of Isabelle is an
intuitionistic fragment of higher-order logic (or simple type theory)
called $\Meta$; its original version was described in \cite{Paulson1989},
and the addition of “sorts” was reported in \cite{Nipkow-LF-91}.

The only predefined type is $\prop$ (“propositions”); new base types
can be postulated when defining objects logics. Types of higher order can be
assembled using the function space constructor $\fun$.

The type of propositions $\prop$ is equipped with a binary operation
${\implies}$ (“meta-implication”) and a universal “meta-quantifier”
$\ALL$, that are used to represent the object
logic rules. As an example, the axiomatization of first-order logic
postulates a type $\tyo$ of booleans, and Modus Ponens
is written as
\begin{equation}\label{eq:modus-ponens}
  \ALL P\,Q.\ \ [P\limp Q] \implies ([P] \implies [Q]).
\end{equation}
The square brackets (which are omitted in Isabelle theories, as well
as the outermost universal quantifiers) represent
an injection from $\tyo$ into $\prop$. %
A consequence of this representation is that every formula of
the object logic appears atomic to $\Meta$.

Types in Isabelle are organized into \emph{classes} and \emph{sorts};
for ease of exposition, we will omit the former.  The axiomatization
of first-order logic postulates a sort $\{\type{term}\}$ (of
“individuals,” or elements of a first-order universe of discourse) and
stipulates that every further type variable $\alpha$ must be of that
sort. In particular, Isabelle/ZF only postulates one new type $\tyi$
(“sets”) of sort $\{\type{term}\}$. Hence, from the type of the universal
quantifier functional $\forall :: (\alpha \fun \tyo) \fun \tyo$, it
follows that it may only be applied to predicates with a variable of
type $\tyi$. This ensures that the object logic is effectively
first-order.

Paulson  \cite{Paulson1989} carried out a proof that the encoding
$\Meta_{\mathrm{IFOL}}$ of
intuitionistic first-order logic IFOL without equality  in the original $\Meta$ is
conservative (there is a correspondence between provable $\phi$ in
IFOL and provable $[\phi]$ in $\Meta_{\mathrm{IFOL}}$) by putting
$\Meta_{\mathrm{IFOL}}$ proofs in \emph{expanded normal form}
\cite{MR0387024}; atomicity as stated after
Equation~(\ref{eq:modus-ponens}) plays a role in this argument. Passing to classical logic does not present
difficulties, but the addition of meta-equality must be taken care of.
Even more so, since the treatment of equality differs between the
original and the present incarnation of $\Meta$; details for the
latter are exhaustively expounded in the recent formalization by
Nipkow and Roßkof \cite{10.1007/978-3-030-79876-5_6}.

The meta-logic $\Meta$ is rather weak; it has no induction/recursion
principles. Types are not inductively presented and, in particular, it
is not possible to prove by induction statements about
object-logic formulas (which are construed as terms of type $\tyi \fun
\dots \fun \tyi \fun \tyo$). Two ways to overcome this limitation are:
\begin{enumerate}
\item
  to construct the
  proof of each instance of the statement by hand or by programming on
  ML; or 
\item
  to encode formulas as sets and prove an internal version statement
  using induction of $\ZF$.
\end{enumerate}

For recursive definitions, only the second option is available, and
that is the way the definition of the forcing relation is implemented
in our formalization.

\subsection{Isabelle/ZF}
\label{sec:isabellezf}

For the most part, the development of set theory in Isabelle is
carried out using its ZF object logic
\cite{DBLP:journals/jar/PaulsonG96}, which is the third logical layer
of the formalization and the most versatile one, since 
Isabelle's native automation is available at this level. Apart from
the type and sort
declarations detailed above, it features a finite axiomatization,
with a predicate for membership, constants for the empty set and an
infinite set, and functions $\isa{Pow}::\tyi\fun\tyi$,
$\union::\tyi\fun\tyi$, and $\isa{PrimReplace} :: \tyi \fun (\tyi
\fun \tyi \fun\tyo)\fun \tyi$ (for Replacement). The Axiom of
Replacement
has a free predicate variable $P$: %
\begin{multline*}
  (\forall x \in A .\ \forall y\, z.\ P(x, y) \wedge P(x, z)
  \longrightarrow y=z) \implies \\
  b \in \isa{PrimReplace}(A, P)
  \longleftrightarrow(\exists x \in A .\ P(x, b)) 
\end{multline*}
The restrictions on sorts described above ensure that it is not
possible that higher-order quantification gets used in $P$. The
statement of $\AC$ also uses a free higher-order variable to denote
an indexed family of nonempty sets. %

The amount of set theory developed in “vanilla” Isabelle/ZF is quite
modest, essentially going no further than Hessenberg's theorem on the
cardinality of products $|A\times A|=|A|\cdot|A| = |A|$. Instead, our decision
(during 2017) to
use this assistant was triggered by its constructibility
library, \session{ZF-Constructible} \citep{paulson_2003},
which contains the development of $L$, the proof that it satisfies
$\AC$, and a version of the Reflection Principle. The latter was
actually encoded as a series of instructions to Isabelle automatic
proof tools that would prove each particular instance of reflection:
This is an example of what was said at the end of Section~\ref{sec:isabelle-metalogic-meta}.

The development of relativization and absoluteness for classes $C::
\isa{i} \fun \isa{o}$ follows the same pattern. Each particular
concept was manually written in a relational form and relativized.
Here, the contrast between the usual way one regards $\ZF$ as a
first-order theory in the language $\{\in \}$ and the mathematical
practice of freely using defined concepts comes to the
forefront. Assistants have refined mechanisms to cope with defined
concepts and the introduction of new notation (which also make their foundations more complicated
than plain first-order logic), and this is the only way that nontrivial
mathematics can be formalized. But when studying relative interpretations, one
usually assumes a spartan syntax and defines relativization by
induction on formulas of the more succinct language. The approach
taken in \session{ZF-Constructible} is to consider relativizations of
the formulas that define each concept. For instance,
in the case of unions, we find a relativization
$\isa{big{\uscore}union}:: (\isa{i} \fun \isa{o}) \fun \isa{i}
\fun \isa{i} \fun \isa{o}$ of the statement
“$\union A = z$”:
\[
 \isa{big{\uscore}union}(M,A,z) \equiv \forall x[M].\ x \in z
 \longleftrightarrow (\exists y[M].\ y \in A \land x \in y)
\]
where $\forall x[M]\dots$ stands for $\forall x.\ M(x)\limp \dots$,
etc. The need to work with \emph{relational} presentations of defined
concepts stems from the fact that the model-theoretic definition of
$L$ requires working with set models and satisfaction, which is
defined for (“codes” of) formulas in the language $\{\in \}$
(viz.\ next Section~\ref{sec:internalized-formulas}).

There is one further point concerning the organization of
\session{ZF-Constructible}. Isabelle provides a very convenient way to
define “contexts,” called \emph{locales}, in which some variables are
fixed and assumptions are made. In the case of the constructibility
library, several locales are defined where the variable $M$ is assumed to
denote a class satisfying certain finite amount of $\ZF$; the weakest
one, \locale{M{\uscore}trans} \cite[Sect.~3]{2020arXiv200109715G}, just
assumes that $M$ is transitive and nonempty. Inside such context, many
absoluteness results are proved. In order to quote those results for a
particular class $C$, one has to \emph{interpret} the locale at
$C$, which amounts to prove that $C$ satisfies the assumptions made by
the context.

\subsection{Internalized formulas}
\label{sec:internalized-formulas}

\session{ZF-Constructible} defines the set $\formula$ of first-order 
formulas in the language $\{ \in \}$, internalized as sets.%
\footnote{These, alongside with lists, are instances of
Isabelle/ZF treatment of inductively defined (internal) datatypes \cite[Sect.~4]{Paulson1995-wz};
induction and recursion theorems for them are proved automatically
(this is in contrast to general well-founded recursion, for which
one has to work with the fundamental recursor $\isa{wfrec}$).}  Its
atomic formulas have the form $\cdot x \in y\cdot$ and $\cdot x =
y\cdot$. (We use dots as a visual aid signaling internalized formulas.)
Variables are represented by de Bruijn indices \cite{MR0321704}, so in
those formulas $x,y \in \omega$; for $\phi\in\formula$ and
$z\in\omega$, $z$ is free in $\phi$ if it occurs under at most $z$
quantifiers. The \isa{arity} function on $\phi$ is one plus the
maximum free index occurring in $\phi$.

The satisfaction predicate
$\isa{sats}::\tyi\fun\tyi\fun\tyi\fun\tyo$ takes as arguments a set
$M$, a list $\mathit{env}\in\isa{list}(M)$ for the assignment of
free indices, and $\phi\in\formula$, and it is written
$M,\mathit{env}\models\phi$ in our formalization.
This completes the
description of the fourth and last formal layer of the development.

Internalized formulas for most (but not all) of the relational
concepts can be obtained by guiding the automatic tactics. But in the
early development of 
\session{ZF-Constructible}, most of the concepts were internalized by
hand; this is the case for union,
\begin{multline*}
  \isa{big{\uscore}union{\uscore}fm}(A,z) \equiv \\
  (({\cdot}\forall{\cdot}\,{\cdot}0 \in \mathit{succ}(z){\cdot} \longleftrightarrow
  ({\cdot}\exists{\cdot}\,{\cdot}0 \in \mathit{succ}(\mathit{succ}(A)){\cdot} \land {\cdot}1 \in
  0{\cdot}\,{\cdot}{\cdot}){\cdot}{\cdot})
\end{multline*}
for which we have the following satisfaction lemma:
\begin{multline}\label{eq:sats_big_union_fm}
  A \in \omega \implies z \in \omega \implies \mathit{env} \in \isa{list}(M)
  \implies \\
  \bigl(M, \mathit{env} \models \isa{big{\uscore}union{\uscore}fm}(A,z)\bigr)%
  \longleftrightarrow \\
  \isa{big{\uscore}union}(\isa{\#\#} M, \isa{nth}(A,
  \mathit{env}), \isa{nth}(z, \mathit{env}))
\end{multline}
Above, $\isa{nth}(x,\mathit{env})$ is the $x$th element of $\mathit{env}$
and $\isa{\#\#}M::\tyi\fun \tyo$ is the class corresponding to the
set $M::\tyi$.

\section{Relative versions of non-absolute concepts}
\label{sec:relat-vers-non-absol}

The treatment of relativization/internalization described in the
previous sections was enough for Paulson's treatment of
constructibility. This is the case because essentially all the
concepts in the way of proving the consistency of $\AC$ are
absolute, and the treatment of relational versions and relativized notions
could be minimized after proving the relevant absoluteness results:
For example, the lemma \isa{Union{\uscore}abs},
\[
  M(A) \implies M(z) \implies \isa{big{\uscore}union}(M, A, z) \longleftrightarrow z = \union
  A
\]
proved under the assumption that $M$ is transitive and nonempty.

Our first attempt to relativize cardinal arithmetic proceeded in the
same way
and we rapidly found out that stating and proving statements like $(||A||
= |A|) ^M$ in a completely relational language was extremely
cumbersome. This observation lead to the discovery of the discipline
expounded in the next subsection.

\subsection{Discipline and tools for relativization}
\label{sec:tools-relativization}
The missing step, that naturally appears in the literature, consists
of having relative \emph{functions} like $\Pow^M$, and the ability to
translate between the different presentations discussed so far.

To achieve this, we provide automatic tools to ease the definitions of
such relative versions, their fully relational counterparts, and the
internalized formulas. For instance, consider the
$\isa{cardinal}::\tyi \fun \tyi$ function defined in
\session{Isabelle/ZF}. Then the commands
\begin{isabelle}
  \isacommand{relativize}\isamarkupfalse%
  \ \isakeyword{functional}\ {\isachardoublequoteopen}cardinal{\isachardoublequoteclose}\ {\isachardoublequoteopen}cardinal{\isacharunderscore}{\kern0pt}rel{\isachardoublequoteclose}\ \isakeyword{external}\isanewline
  \isacommand{relationalize}\isamarkupfalse%
  \ {\isachardoublequoteopen}cardinal{\isacharunderscore}{\kern0pt}rel{\isachardoublequoteclose}\ {\isachardoublequoteopen}is{\isacharunderscore}{\kern0pt}cardinal{\isachardoublequoteclose}\isanewline
  \isacommand{synthesize}\isamarkupfalse%
  \ {\isachardoublequoteopen}is{\isacharunderscore}{\kern0pt}cardinal{\isachardoublequoteclose}\ \isakeyword{from{\isacharunderscore}{\kern0pt}definition}\ \isakeyword{assuming}\ {\isachardoublequoteopen}nonempty{\isachardoublequoteclose}%
\end{isabelle}
define the relative cardinal function
$\isa{cardinal{\uscore}rel}::(\tyi \fun \tyo) \fun \tyi \fun\tyi$
(whose first argument is the class $M$, and is denoted $|\cdot|^M$ as expected),
the relational version $\isa{is{\uscore}cardinal}$ of the latter, the
internalized formula \isa{is{\uscore}cardinal{\uscore}fm} whose
satisfaction by a set is equivalent to the relational version, and
prove the previous statement (analogous to (\ref{eq:sats_big_union_fm})).
The proof that $\isa{is{\uscore}cardinal}(M,x,z)$  encodes the
statement $|x|^M = z$ must still be done by hand, since the definition
of $\isa{cardinal{\uscore}rel}$ already involves some tacit
absoluteness results (“\textit{the least $z \in \Ord$ such that $z
  \approx^M x$}” instead
of “\textit{the least $z \in \Ord^M$ such that $z
  \approx^M x$}”, and the like).

\subsection{Extension of Isabelle/ZF}
\label{sec:extension-isabellezf}
We extended \cite{Delta_System_Lemma-AFP} the material formalized in
Isabelle, from basic results involving function spaces and the
definition of cardinal exponentiation, to a treatment of cofinality
and the Delta System Lemma for $\omega_1$-families. We also included a
concise treatment of the axiom of Dependent Choices $\DC$ and the
general version of Rasiowa-Sikorski Lemma \cite{2018arXiv180705174G}
and a choiceless one for countable preorders.

This material was subsequently put in relative form in our formal
development on transitive class models \cite{Transitive_Models-AFP}
using as an aid the tools from
Section~\ref{sec:tools-relativization}. We also relativized many
original theories appearing in Isabelle/ZF, including the
fundamentals of cardinal arithmetic, the cumulative hierarchy, and the
definition of the $\ale{}$ function.

\section{Set models and forcing}
\label{sec:forcing}

\subsection{The $\ZFC$ axioms as locales}\label{sec:zfc-axioms-as-locales}
The description of set models of fragments of $\ZFC$ was performed
using locales that fix a variable $M::\tyi$ and pack assumptions
stating that $\lb M, \in\rb$ satisfy some axioms; for example, the
locale \locale{M{\uscore}Z{\uscore}basic} states that Zermelo set
theory holds in $M$. It would be natural to state those assumptions
directly as the corresponding satisfactions, as in
\[
  M, [\,]\models \isa{{\isasymcdot}Union\ Ax{\isasymcdot}}
\]
where \isa{{\isasymcdot}Union\ Ax{\isasymcdot}} is the $\formula$ code for the
Union Axiom. We actually decided to express the axioms other than the
infinite schemes in relational form, by using terms already
available in \session{ZF-Constructible} and for which useful lemmas
had already been proved (and, as it was mentioned in Section~\ref{sec:isabellezf},
this third layer of the formalization has more tools at our
disposal); the Union Axiom (\isa{Union{\uscore}ax}), for
instance, is defined as follows:
\[
\forall x[\isa{\#\#}M].\ \exists z[\isa{\#\#}M].\ \isa{big{\uscore}union}(\isa{\#\#}M,x,z)
\]
Both assumptions are then shown to be equivalent:
\[
  \isa{Union{\uscore}ax}(\isa{\#\#}M) \lsii M, [\,]\models \isa{{\isasymcdot}Union\ Ax{\isasymcdot}}
\]

For stating the axiom schemes, \session{ZF-Constructible} defines
$\isa{separation}(N,Q)$:
\[
  \forall z[N].\ \exists y[N].\ \forall
  x[N].\ x\in y \longleftrightarrow x\in z \land Q(x),
\]
and $\isa{strong{\uscore}replacement}(N,R)$:
\begin{multline*}
  \forall
  A[N].\ \isa{univalent}(N,A,R) \longrightarrow \\
  (\exists Y[N].\ \forall b[N].\ b\in Y\longleftrightarrow (\exists x[N].\ x\in
  A \land R(x,b))),
\end{multline*}
whose first argument $N$ is a class and their second arguments $Q$ and
$R$ are predicates of types $Q :: \tyi \fun \tyo$ and $R ::
\tyi \fun \tyi \fun \tyo$, respectively (the \isatt{univalent}
predicate states that $R$ is functional in its first variable; see Appendix~\ref{sec:def-main-relative}). The Separation Axiom appears in
\locale{M{\uscore}Z{\uscore}basic} as follows, where 
\isa{\isasymphi} is free and \isa{\isacharat} denotes
list concatenation:
\begin{isabelle}
separation{\isacharunderscore}{\kern0pt}ax{\isacharcolon}{\kern0pt}\ {\isachardoublequoteopen}{\isasymphi}\ {\isasymin}\ formula\ {\isasymLongrightarrow}\ env\ {\isasymin}\ list{\isacharparenleft}{\kern0pt}M{\isacharparenright}{\kern0pt}\ {\isasymLongrightarrow}\isanewline
\ \ \ \ \ \ \ \ \ \ \ \ \ \ \ \ \ \ \ \ arity{\isacharparenleft}{\kern0pt}{\isasymphi}{\isacharparenright}{\kern0pt}\ {\isasymle}\ {\isadigit{1}}\ {\isacharplus}{\kern0pt}\isactrlsub {\isasymomega}\ length{\isacharparenleft}{\kern0pt}env{\isacharparenright}{\kern0pt}\ {\isasymLongrightarrow}\isanewline
\ \ \ \ \ \ \ \ \ \ \ \ \ \ \ \ \ \ \ \ separation{\isacharparenleft}{\kern0pt}{\isacharhash}{\kern0pt}{\isacharhash}{\kern0pt}M{\isacharcomma}{\kern0pt}{\isasymlambda}x{\isachardot}{\kern0pt}\ {\isacharparenleft}{\kern0pt}M{\isacharcomma}{\kern0pt}\ {\isacharbrackleft}{\kern0pt}x{\isacharbrackright}{\kern0pt}\ {\isacharat}{\kern0pt}\ env\ {\isasymTurnstile}\ {\isasymphi}{\isacharparenright}{\kern0pt}{\isacharparenright}{\kern0pt}{\isachardoublequoteclose}
\end{isabelle}
Note that the predicate $Q$ mentioned above corresponds to the
satisfaction of \isa{\isasymphi}.

In contrast to Separation, we stated each instance of the Replacement
Axiom separately by means of the
$\isa{replacement{\uscore}assm}(M,\mathit{env},\phi)$ predicate:
\begin{isabelle}
{\isasymphi}\ {\isasymin}\ formula\ {\isasymlongrightarrow}\ env\ {\isasymin}\ list{\isacharparenleft}{\kern0pt}M{\isacharparenright}{\kern0pt}\ {\isasymlongrightarrow}\isanewline
\ \ arity{\isacharparenleft}{\kern0pt}{\isasymphi}{\isacharparenright}{\kern0pt}\ {\isasymle}\ {\isadigit{2}}\ {\isacharplus}{\kern0pt}\isactrlsub {\isasymomega}\ length{\isacharparenleft}{\kern0pt}env{\isacharparenright}{\kern0pt}\ {\isasymlongrightarrow}\isanewline
\ \ \ \ strong{\isacharunderscore}{\kern0pt}replacement{\isacharparenleft}{\kern0pt}{\isacharhash}{\kern0pt}{\isacharhash}{\kern0pt}M{\isacharcomma}{\kern0pt}{\isasymlambda}x\ y{\isachardot}{\kern0pt}\ {\isacharparenleft}{\kern0pt}M\ {\isacharcomma}{\kern0pt}\ {\isacharbrackleft}{\kern0pt}x{\isacharcomma}{\kern0pt}y{\isacharbrackright}{\kern0pt}{\isacharat}{\kern0pt}env\ {\isasymTurnstile}\ {\isasymphi}{\isacharparenright}{\kern0pt}{\isacharparenright}{\kern0pt}{\isachardoublequoteclose},
\end{isabelle}
which encodes the instance of replacement for the model $M$
corresponding to the predicate $R(x,y)$ given by
$\phi^M(x,y,x_0,\dots,x_n)$, where $\isatt{env} = [x_0,\dots,x_n]$.

In turn, the \isa{{\isasymcdot}Replacement{\isasymcdot}} function
takes a formula code and returns the corresponding replacement
instance:
\begin{isabelle}
{\isasymphi}\ {\isasymin}\ formula\ {\isasymLongrightarrow}\isanewline
{\isacharparenleft}{\kern0pt}M{\isacharcomma}{\kern0pt}\ {\isacharbrackleft}{\kern0pt}{\isacharbrackright}{\kern0pt}\ {\isasymTurnstile}\ {\isasymcdot}Replacement{\isacharparenleft}{\kern0pt}{\isasymphi}{\isacharparenright}{\kern0pt}{\isasymcdot}{\isacharparenright}{\kern0pt}\ {\isasymlongleftrightarrow}\ {\isacharparenleft}{\kern0pt}{\isasymforall}env{\isachardot}{\kern0pt}\ replacement{\isacharunderscore}{\kern0pt}assm{\isacharparenleft}{\kern0pt}M{\isacharcomma}{\kern0pt}\ env{\isacharcomma}{\kern0pt}\ {\isasymphi}{\isacharparenright}{\kern0pt}{\isacharparenright}
\end{isabelle}

Starting from \locale{M{\uscore}Z{\uscore}basic}, stronger locales are
defined by assuming more replacement instances.
These assumptions are then invoked to interpret at the class
$\isa{\#\#} M$ the relevant locales appearing in
\session{ZF-Constructible}, and further ones required for the relative
results from Section~\ref{sec:extension-isabellezf}. See
Section~\ref{sec:repl-instances} for details.

\subsection{The fundamental theorems}\label{sec:fundamental-theorems}
At this point, we work inside the locale \isa{M{\uscore}ctm1} that
assumes $M$ to be
countable and transitive, and satisfies some fragment of
$\ZFC$%
\footnote{%
Namely, Zermelo set theory plus the 7 replacement instances included in
the locales \locale{M{\uscore}ZF1} and \locale{M{\uscore}ZF{\uscore}ground}.
}.
This is further extended by assuming a forcing notion
$\lb \PP, {\preceq} ,\1\rb \in M$. The actual implementation reads:

\begin{isabelle}
\isacommand{locale}\isamarkupfalse%
\ forcing{\isacharunderscore}{\kern0pt}notion\ {\isacharequal}{\kern0pt}\isanewline
\ \ \isakeyword{fixes}\ P\ {\isacharparenleft}{\kern0pt}{\isacartoucheopen}{\isasymbbbP}{\isacartoucheclose}{\isacharparenright}{\kern0pt}\ \isakeyword{and}\ leq\ \isakeyword{and}\ one\ {\isacharparenleft}{\kern0pt}{\isacartoucheopen}{\isasymone}{\isacartoucheclose}{\isacharparenright}{\kern0pt}\isanewline
\ \ \isakeyword{assumes}\ one{\isacharunderscore}{\kern0pt}in{\isacharunderscore}{\kern0pt}P{\isacharcolon}{\kern0pt}\ \ \ \ \ \ \ {\isachardoublequoteopen}{\isasymone}\ {\isasymin}\ {\isasymbbbP}{\isachardoublequoteclose}\isanewline
\ \ \ \ \isakeyword{and}\ leq{\isacharunderscore}{\kern0pt}preord{\isacharcolon}{\kern0pt}\ \ \ \ \ \ \ {\isachardoublequoteopen}preorder{\isacharunderscore}{\kern0pt}on{\isacharparenleft}{\kern0pt}{\isasymbbbP}{\isacharcomma}{\kern0pt}leq{\isacharparenright}{\kern0pt}{\isachardoublequoteclose}\isanewline
\ \ \ \ \isakeyword{and}\ one{\isacharunderscore}{\kern0pt}max{\isacharcolon}{\kern0pt}\ \ \ \ \ \ \ \ \ \ {\isachardoublequoteopen}{\isasymforall}p{\isasymin}{\isasymbbbP}{\isachardot}{\kern0pt}\ {\isasymlangle}p{\isacharcomma}{\kern0pt}{\isasymone}{\isasymrangle}{\isasymin}leq{\isachardoublequoteclose}
\end{isabelle}

\begin{isabelle}
\isacommand{locale}\isamarkupfalse%
\ forcing{\isacharunderscore}{\kern0pt}data{\isadigit{1}}\ {\isacharequal}{\kern0pt}\ forcing{\isacharunderscore}{\kern0pt}notion\ {\isacharplus}{\kern0pt}\ M{\isacharunderscore}{\kern0pt}ctm{\isadigit{1}}\ {\isacharplus}{\kern0pt}\isanewline
\ \ \isakeyword{assumes}\ P{\isacharunderscore}{\kern0pt}in{\isacharunderscore}{\kern0pt}M{\isacharcolon}{\kern0pt}\ \ \ \ \ \ \ \ \ \ \ {\isachardoublequoteopen}{\isasymbbbP}\ {\isasymin}\ M{\isachardoublequoteclose}\isanewline
\ \ \ \ \isakeyword{and}\ leq{\isacharunderscore}{\kern0pt}in{\isacharunderscore}{\kern0pt}M{\isacharcolon}{\kern0pt}\ \ \ \ \ \ \ \ \ {\isachardoublequoteopen}leq\ {\isasymin}\ M{\isachardoublequoteclose}
\end{isabelle}
The version of the Forcing Theorems that we formalized follows the
considerations on the $\forces^*$ relation as discussed in Kunen's new
\emph{Set Theory}
\cite[p.~257ff]{kunen2011set}.
We defined forcing for atomic formulas by recursion on names in an
analogous fashion. But, in contrast to the point made on
p.~260 of this book, the structural recursion used to define the forcing
relation was replaced by one  involving codes for formulas. Thus, the metatheoretic formula
transformer $\phi\mapsto\mathit{Forces}_\phi$ was replaced by a
set-theoretic class function $\forceisa:: \tyi \fun \tyi$, which was defined by using
Isabelle/ZF facilities for primitive recursion.

Next, we state this version of the fundamental theorems in a compact
way. For any $G\subseteq \PP$, our notation for the extension of $M$ by
$G$ is the  customary one: $M[G]\defi
\{ \val(G,\tau) : \tau\in M \}$, where the interpretation
$\val(G,\tau)$ is defined by well-founded recursion on
$\tau$.

\begin{theorem}\label{th:forcing-thms}
  For every
  $\phi\in\formula$ with $\isa{arity}(\phi)\leq n$ and $\tau_1,\dots,\tau_n\in M$,
  \begin{enumerate}
  \item\label{item:definability} (Definability)
    $\forceisa(\phi)\in\formula$, 
  \end{enumerate}
  where the 
  arity of $\forceisa(\phi)$ is at most $\isa{arity}(\phi) + 4$; and if
  “$p \forces \phi\ [\tau_1,\dots,\tau_n]$”
  denotes
  “$M, [p,\PP,\preceq,\1, \tau_1,\dots,\tau_n]  \models
  \forceisa(\phi)$”, we have:
  \begin{enumerate}
    \setcounter{enumi}{1}
  \item\label{item:truth-lemma} (Truth Lemma) for every $M$-generic $G$,
    \[
      \exists p\in G.\ \; p \forces \phi\ [\tau_1,\dots,\tau_n]
    \]
    is equivalent to 
    \[
      M[G], [\val(G,\tau_1),\dots,\val(G,\tau_n)]
      \models \phi.
    \]
  \item \label{item:density-lemma} (Density Lemma) $p \forces \phi\ [\tau_1,\dots,\tau_n]$
    if and only if 
    $\{q\in \PP :  q \forces \phi\ [\tau_1,\dots,\tau_n]\}$
    is dense below $p$.
  \end{enumerate}
\end{theorem}
The items in Theorem~\ref{th:forcing-thms} appear in our
\session{Independence\_CH} session \cite{Independence_CH-AFP} as three
separate lemmas (located in the theory
\theory{Forcing{\uscore}Theorems}).
For instance, the Truth Lemma is stated as
follows:
\begin{isabelle}
\isacommand{lemma}\isamarkupfalse%
\ truth{\isacharunderscore}{\kern0pt}lemma{\isacharcolon}{\kern0pt}\isanewline
\ \ \isakeyword{assumes}\isanewline
\ \ \ \ {\isachardoublequoteopen}{\isasymphi}{\isasymin}formula{\isachardoublequoteclose}\isanewline
\ \ \ \ {\isachardoublequoteopen}env{\isasymin}list{\isacharparenleft}{\kern0pt}M{\isacharparenright}{\kern0pt}{\isachardoublequoteclose}\ {\isachardoublequoteopen}arity{\isacharparenleft}{\kern0pt}{\isasymphi}{\isacharparenright}{\kern0pt}{\isasymle}length{\isacharparenleft}{\kern0pt}env{\isacharparenright}{\kern0pt}{\isachardoublequoteclose}\isanewline
\ \ \isakeyword{shows}\isanewline
\ \ \ \ {\isachardoublequoteopen}{\isacharparenleft}{\kern0pt}{\isasymexists}p{\isasymin}G{\isachardot}{\kern0pt}\ p\ {\isasymtturnstile}\ {\isasymphi}\ env{\isacharparenright}{\kern0pt}\ \ \ {\isasymlongleftrightarrow}\ \ \ M{\isacharbrackleft}{\kern0pt}G{\isacharbrackright}{\kern0pt}{\isacharcomma}{\kern0pt}\ map{\isacharparenleft}{\kern0pt}val{\isacharparenleft}{\kern0pt}G{\isacharparenright}{\kern0pt}{\isacharcomma}{\kern0pt}env{\isacharparenright}{\kern0pt}\ {\isasymTurnstile}\ {\isasymphi}{\isachardoublequoteclose}
\end{isabelle}
where the $\forces$ notation (and its precedence) had already been set up in the
\theory{Forces{\uscore}Definition} theory.

Kunen first describes forcing for atomic formulas using a mutual
recursion
but then \cite[p.~257]{kunen2011set} it is cast as a single
recursively defined function $F$ over a well-founded  relation $R$.
In our formalization, these are called $\frcat$ and 
$\isa{frecR}$, respectively, and are defined on tuples $\lb \mathit{ft},t_1,t_2,p\rb$ (where
$\mathit{ft}\in\{0,1\}$ indicates whether the atomic formula being
forced is an equality or a membership, respectively).
Forcing for general formulas is then defined by recursion on the
datatype $\formula$ as indicated above. Technical details on the
implementation and proofs of the
Forcing Theorems have been spelled out in our
\cite{2020arXiv200109715G}.

\section{A sample formal proof}
\label{sec:sample-formal-proof}

We present a fragment of the formal version of the proof that the
Powerset Axiom holds in a generic extension, which also serves to
illustrate the Isar dialect of Isabelle.

We quote the relevant
paragraph of Kunen's \cite[Thm.~IV.2.27]{kunen2011set}:
\begin{quote}
  For Power Set (similarly to Union above), it is sufficient to prove
  that whenever $a \in M[G]$, there is a $b \in M[G]$ such that
  $\mathcal{P}(a) \cap M[G] \subseteq b$. Fix $\tau \in
  M^{\mathbb{P}}$ such that $\tau_{G}=a$. Let
  $Q=(\mathcal{P}(\operatorname{dom}(\tau) \times
  \mathbb{P}))^{M}$. This is the set of all names $\vartheta \in
  M^{\mathbb{P}}$ such that $\operatorname{dom}(\vartheta) \subseteq
  \operatorname{dom}(\tau)$. Let $\pi=Q \times\{\1\}$ and let
  $b=\pi_{G}=$ $\left\{\vartheta_{G}: \vartheta \in Q\right\}$. Now,
  consider any $c \in \mathcal{P}(a) \cap M[G]$; we need to show that
  $c \in b$. Fix $\chi \in M^{\mathbb{P}}$ such that
  $\chi_{G}=c$, and let $\vartheta=\{\langle\sigma, p\rangle:
  \sigma \in \operatorname{dom}(\tau) \wedge p \Vdash \sigma \in
  \chi\}$; $\vartheta \in M$ by the Definability Lemma. Since
  $\vartheta \in Q$, we are done if we can show that
  $\vartheta_{G}=c$.
\end{quote}
The assumption $a\in M[G]$ appears in the lemma statement, and the
goal involving $b$ in the first sentence will appear below (signaled
by “{\small (**)}”); formalized
material necessarily tends to be much more linear than usual prose. In
what follows, we
will intersperse the relevant passages of the proof.
\begin{isabelle}
\isacommand{lemma}\isamarkupfalse%
\ Pow{\isacharunderscore}{\kern0pt}inter{\isacharunderscore}{\kern0pt}MG{\isacharcolon}{\kern0pt}\isanewline
\ \ \isakeyword{assumes}\ {\isachardoublequoteopen}a{\isasymin}M{\isacharbrackleft}{\kern0pt}G{\isacharbrackright}{\kern0pt}{\isachardoublequoteclose}\isanewline
\ \ \isakeyword{shows}\ {\isachardoublequoteopen}Pow{\isacharparenleft}{\kern0pt}a{\isacharparenright}{\kern0pt}\ {\isasyminter}\ M{\isacharbrackleft}{\kern0pt}G{\isacharbrackright}{\kern0pt}\ {\isasymin}\ M{\isacharbrackleft}{\kern0pt}G{\isacharbrackright}{\kern0pt}{\isachardoublequoteclose}\isanewline
\isacommand{proof}\isamarkupfalse%
\ {\isacharminus}{\kern0pt}
\end{isabelle}
\textit{Fix $\tau \in  M^{\mathbb{P}}$ such that $\tau_{G}=a$.}
\begin{isabelle}
\ \ \isacommand{from}\isamarkupfalse%
\ assms\isanewline
\ \ \isacommand{obtain}\isamarkupfalse%
\ {\isasymtau}\ \isakeyword{where}\ {\isachardoublequoteopen}{\isasymtau}\ {\isasymin}\ M{\isachardoublequoteclose}\ {\isachardoublequoteopen}val{\isacharparenleft}{\kern0pt}G{\isacharcomma}{\kern0pt}\ {\isasymtau}{\isacharparenright}{\kern0pt}\ {\isacharequal}{\kern0pt}\ a{\isachardoublequoteclose}\isanewline
\ \ \ \ \isacommand{using}\isamarkupfalse%
\ GenExtD\ \isacommand{by}\isamarkupfalse%
\ auto
\end{isabelle}
\textit{Let
  $Q=(\mathcal{P}(\operatorname{dom}(\tau) \times
  \mathbb{P}))^{M}$. This is the set of all names $\vartheta \in
  M^{\mathbb{P}}$} [\dots]%
\begin{isabelle}
\ \ \isacommand{let}\isamarkupfalse%
\ {\isacharquery}{\kern0pt}Q{\isacharequal}{\kern0pt}{\isachardoublequoteopen}Pow\isactrlbsup M\isactrlesup {\isacharparenleft}{\kern0pt}domain{\isacharparenleft}{\kern0pt}{\isasymtau}{\isacharparenright}{\kern0pt}{\isasymtimes}{\isasymbbbP}{\isacharparenright}{\kern0pt}{\isachardoublequoteclose}
\end{isabelle}
\textit{ Let $\pi=Q \times\{\1\}$ and let
  $b=\pi_{G}=$ $\left\{\vartheta_{G}: \vartheta \in Q\right\}$.}
\begin{isabelle}
\ \ \isacommand{let}\isamarkupfalse%
\ {\isacharquery}{\kern0pt}{\isasympi}{\isacharequal}{\kern0pt}{\isachardoublequoteopen}{\isacharquery}{\kern0pt}Q{\isasymtimes}{\isacharbraceleft}{\kern0pt}{\isasymone}{\isacharbraceright}{\kern0pt}{\isachardoublequoteclose}\isanewline
\ \ \isacommand{let}\isamarkupfalse%
\ {\isacharquery}{\kern0pt}b{\isacharequal}{\kern0pt}{\isachardoublequoteopen}val{\isacharparenleft}{\kern0pt}G{\isacharcomma}{\kern0pt}{\isacharquery}{\kern0pt}{\isasympi}{\isacharparenright}{\kern0pt}{\isachardoublequoteclose}
\end{isabelle}
(Recall: \textit{\dots there is a $b\in M[G]$ such that\dots})
\begin{isabelle}
\ \ \isacommand{from}\isamarkupfalse%
\ {\isacartoucheopen}{\isasymtau}{\isasymin}M{\isacartoucheclose}\isanewline
\ \ \isacommand{have}\isamarkupfalse%
\ {\isachardoublequoteopen}domain{\isacharparenleft}{\kern0pt}{\isasymtau}{\isacharparenright}{\kern0pt}{\isasymtimes}{\isasymbbbP}\ {\isasymin}\ M{\isachardoublequoteclose}\ {\isachardoublequoteopen}domain{\isacharparenleft}{\kern0pt}{\isasymtau}{\isacharparenright}{\kern0pt}\ {\isasymin}\ M{\isachardoublequoteclose}\isanewline
\ \ \ \ \isacommand{by}\isamarkupfalse%
\ simp{\isacharunderscore}{\kern0pt}all\isanewline
\ \ \isacommand{then}\isamarkupfalse%
\isanewline
\ \ \isacommand{have}\isamarkupfalse%
\ {\isachardoublequoteopen}{\isacharquery}{\kern0pt}b\ {\isasymin}\ M{\isacharbrackleft}{\kern0pt}G{\isacharbrackright}{\kern0pt}{\isachardoublequoteclose}\isanewline
\ \ \ \ \isacommand{by}\isamarkupfalse%
\ {\isacharparenleft}{\kern0pt}auto\ intro{\isacharbang}{\kern0pt}{\isacharcolon}{\kern0pt}GenExtI{\isacharparenright}{\kern0pt}
\end{isabelle}
\textit{Now,
  consider any $c \in \mathcal{P}(a) \cap M[G]$; we need to show that
  $c \in b$.}
\begin{isabelle}
  \label{goal-on-b}
\ \ \isacommand{have}\isamarkupfalse%
\ {\isachardoublequoteopen}Pow{\isacharparenleft}{\kern0pt}a{\isacharparenright}{\kern0pt}\ {\isasyminter}\ M{\isacharbrackleft}{\kern0pt}G{\isacharbrackright}{\kern0pt}\ {\isasymsubseteq}\ {\isacharquery}{\kern0pt}b{\isachardoublequoteclose}\hfill
\mbox{\rm\small(**)}\isanewline
\ \ \isacommand{proof}\isamarkupfalse%
\isanewline
\ \ \ \ \isacommand{fix}\isamarkupfalse%
\ c\isanewline
\ \ \ \ \isacommand{assume}\isamarkupfalse%
\ {\isachardoublequoteopen}c\ {\isasymin}\ Pow{\isacharparenleft}{\kern0pt}a{\isacharparenright}{\kern0pt}\ {\isasyminter}\ M{\isacharbrackleft}{\kern0pt}G{\isacharbrackright}{\kern0pt}{\isachardoublequoteclose}
\end{isabelle}
\textit{Fix $\chi \in M^{\mathbb{P}}$ such that
  $\chi_{G}=c$,}
\begin{isabelle}
\ \ \ \ \isacommand{then}\isamarkupfalse%
\isanewline
\ \ \ \ \isacommand{obtain}\isamarkupfalse%
\ {\isasymchi}\ \isakeyword{where}\ {\isachardoublequoteopen}c{\isasymin}M{\isacharbrackleft}{\kern0pt}G{\isacharbrackright}{\kern0pt}{\isachardoublequoteclose}\ {\isachardoublequoteopen}{\isasymchi}\ {\isasymin}\ M{\isachardoublequoteclose}\ {\isachardoublequoteopen}val{\isacharparenleft}{\kern0pt}G{\isacharcomma}{\kern0pt}{\isasymchi}{\isacharparenright}{\kern0pt}\ {\isacharequal}{\kern0pt}\ c{\isachardoublequoteclose}\isanewline
\ \ \ \ \ \ \isacommand{using}\isamarkupfalse%
\ GenExt{\isacharunderscore}{\kern0pt}iff\ \isacommand{by}\isamarkupfalse%
\ auto
\end{isabelle}
\textit{and let $\vartheta=\{\langle\sigma, p\rangle:
  \sigma \in \operatorname{dom}(\tau) \wedge p \Vdash \sigma \in
  \chi\}$;}
\begin{isabelle}
\ \ \ \ \isacommand{let}\isamarkupfalse%
\ {\isacharquery}{\kern0pt}{\isasymtheta}{\isacharequal}{\kern0pt}{\isachardoublequoteopen}{\isacharbraceleft}{\kern0pt}{\isasymlangle}{\isasymsigma}{\isacharcomma}{\kern0pt}p{\isasymrangle}\ {\isasymin}domain{\isacharparenleft}{\kern0pt}{\isasymtau}{\isacharparenright}{\kern0pt}{\isasymtimes}\isasymbbbP\ {\isachardot}{\kern0pt}\ p\ {\isasymtturnstile}\ {\isasymcdot}{\isadigit{0}}\ {\isasymin}\ {\isadigit{1}}{\isasymcdot}\ {\isacharbrackleft}{\kern0pt}{\isasymsigma}{\isacharcomma}{\kern0pt}{\isasymchi}{\isacharbrackright}{\kern0pt}\ {\isacharbraceright}{\kern0pt}{\isachardoublequoteclose}
\end{isabelle}
\textit{$\vartheta \in M$ by the Definability Lemma.}
\begin{isabelle}
\ \ \ \ \isacommand{have}\isamarkupfalse%
\ {\isachardoublequoteopen}arity{\isacharparenleft}{\kern0pt}forces{\isacharparenleft}{\kern0pt}\ {\isasymcdot}{\isadigit{0}}\ {\isasymin}\ {\isadigit{1}}{\isasymcdot}\ {\isacharparenright}{\kern0pt}{\isacharparenright}{\kern0pt}\ {\isacharequal}{\kern0pt}\ {\isadigit{6}}{\isachardoublequoteclose}\isanewline
\ \ \ \ \ \ \isacommand{using}\isamarkupfalse%
\ arity{\isacharunderscore}{\kern0pt}forces{\isacharunderscore}{\kern0pt}at\ \isacommand{by}\isamarkupfalse%
\ auto\isanewline
\ \ \ \ \isacommand{with}\isamarkupfalse%
\ {\isacartoucheopen}domain{\isacharparenleft}{\kern0pt}{\isasymtau}{\isacharparenright}{\kern0pt}\ {\isasymin}\ M{\isacartoucheclose}\ {\isacartoucheopen}{\isasymchi}\ {\isasymin}\ M{\isacartoucheclose}\isanewline
\ \ \ \ \isacommand{have}\isamarkupfalse%
\ {\isachardoublequoteopen}{\isacharquery}{\kern0pt}{\isasymtheta}\ {\isasymin}\ M{\isachardoublequoteclose}\isanewline
\ \ \ \ \ \ \isacommand{using}\isamarkupfalse%
\ sats{\isacharunderscore}{\kern0pt}fst{\isacharunderscore}{\kern0pt}snd{\isacharunderscore}{\kern0pt}in{\isacharunderscore}{\kern0pt}M\isanewline
\ \ \ \ \ \ \isacommand{by}\isamarkupfalse%
\ simp\end{isabelle}
\textit{Since
  $\vartheta \in Q$,}
\begin{isabelle}
\ \ \ \ \isacommand{with}\isamarkupfalse%
\ {\isacartoucheopen}domain{\isacharparenleft}{\kern0pt}{\isasymtau}{\isacharparenright}{\kern0pt}{\isasymtimes}{\isasymbbbP}\ {\isasymin}\ M{\isacartoucheclose}\isanewline
\ \ \ \ \isacommand{have}\isamarkupfalse%
\ {\isachardoublequoteopen}{\isacharquery}{\kern0pt}{\isasymtheta}\ {\isasymin}\ {\isacharquery}{\kern0pt}Q{\isachardoublequoteclose}\isanewline
\ \ \ \ \ \ \isacommand{using}\isamarkupfalse%
\ Pow{\isacharunderscore}{\kern0pt}rel{\isacharunderscore}{\kern0pt}char\ \isacommand{by}\isamarkupfalse%
\ auto
\end{isabelle}
\textit{we are done if we can show that
  $\vartheta_{G}=c$.}
\begin{isabelle}
\ \ \ \ \isacommand{have}\isamarkupfalse%
\ {\isachardoublequoteopen}val{\isacharparenleft}{\kern0pt}G{\isacharcomma}{\kern0pt}{\isacharquery}{\kern0pt}{\isasymtheta}{\isacharparenright}{\kern0pt}\ {\isacharequal}{\kern0pt}\ c{\isachardoublequoteclose}\isanewline
\ \ \ \ \isacommand{proof} \ \mbox{$[\dots]$}
\end{isabelle}

This cherry-picked example shows that the formalization can be close
to the mathematical exposition and might be useful to reconstruct the
proof from the book; nonetheless, it also has significantly more
details than the mathematical prose, even with some indications to
direct the automatic tools.

There has been some progress on assistants where one writes statements
and proofs in natural language; recently P.~Koepke and his team
achieved magnificent results by using Isabelle/Naproche
\cite{10.1007/978-3-030-81097-9_2} to formalize proofs of several
results (particularly, the proof of König's Theorem). The input
language of Isabelle/Naproche is a \emph{controlled} natural language
that presents the result being formalized as a deduction in
first-order logic, where every assumption and the ``whole logical
scenario'' are explicitly given. From the input language,
Isabelle/Naproche builds ``proof tasks'' that are handled to automatic
theorem provers. As far as we can tell, Isabelle/Naproche is
promising but still unsuitable for a project of the magnitude of ours.

\section{Main achievements of the formalization}
\label{sec:main-achievements}

\subsection{A sufficient set of replacement instances}
\label{sec:repl-instances}

We isolated 22 instances of Replacement that are sufficient to force
$\CH$ or $\neg\CH$. Many of these were already present in
relational form in the \session{ZF-Constructible} library.

The first 4 instances, collected in the subset
\isa{instances1{\uscore}fms} of \formula, consist of basic
constructions:

\begin{itemize}
\item 2 instances for transitive closure: one to prove closure under
  iteration of $X\mapsto\union X$ and an auxiliary one used to show absoluteness.
\item 1 instance to define $\in$-rank.
\item 1 instance to construct the cumulative hierarchy (rank initial segments).
\end{itemize}

The next 4 instances (gathered in \isa{instances2{\uscore}fms})
are needed to set up
cardinal arithmetic in $M$:
\begin{itemize}
\item 2 instances for the definition of
  ordertypes: The relevant well-founded recursion and a technical
  auxiliary instance.
\item 2 instances for Aleph: Replacement through the ordertype function (for Hartogs' Theorem) and the well-founded recursion
  using it.
\end{itemize}

We also need a one extra replacement instance $\psi$ on $M$ for each
$\phi$ of the
previous ones to have them in $M[G]$:
\[
  \psi(x,\alpha,y_1,\dots,y_n) \defi \quine{\alpha = \min \bigl\{
    \beta \mid \exists\tau\in V_\beta.\  \mathit{snd}(x) \forces
    \phi\ [\mathit{fst}(x),\tau,y_1,\dots,y_n]\bigr\}}
\]
Here, $\mathit{fst}(\lb a,b\rb) = a$ and $\mathit{snd}(\lb a,b\rb) = b$.
The map $\phi\mapsto\psi$ is
the function $\calF$ referred to in the abstract.
All such “ground” replacement
instances appear in the locale \locale{M{\uscore}ZF3} and form the set
\isa{instances3{\uscore}fms}.

That makes 16 instances up to now. For the setup of forcing, we
require the following 3 instances, which form the set
\isa{instances{\uscore}ground{\uscore}fms}:
\begin{itemize}
\item Well-founded recursion to define check-names.
\item Well-founded recursion for the definition of forcing for atomic formulas.
\item Replacement through $x\mapsto \lb x,\check{x}\rb$ (for the
  definition of $\punto{G}$).
\end{itemize}
The proof of the $\Delta$-System Lemma requires 2 instances which form the set
\isa{instances{\uscore}ground{\uscore}notCH{\uscore}fms}, that are
used for the recursive construction of sets using a choice function (as in the
construction of a wellorder of $X$ given a choice function on
$\Pow(X)$), and to show its absoluteness.

The $21$ formulas up to this point are collected into the set
\isa{overhead{\uscore}notCH} (called $\Omega$ in the abstract), which is enough to
force $\neg\CH$. To force $\CH$, we required one further instance for
the absoluteness of the recursive construction in the proof of
Dependent Choices from $\AC$. A listing with the names of all the formulas
can be found in Appendix~\ref{sec:repl-instances-appendix}.
  
The particular choice of some of the instances above arose from
Paulson's architecture on which we based our development.
This applies every time
a locale from \session{ZF-Constructible} has to be
interpreted (\locale{M{\uscore}eclose} and
\locale{M{\uscore}ordertype}, respectively, for the “auxiliary” instances).

On the other hand, we replaced the original proof of the
Schröder-Bernstein Theorem by Zermelo's one
\cite[Exr. x4.27]{moschovakis1994notes}, because the former required
at least one extra instance
arising from an iteration. We also managed to avoid 12 further
replacements by restructuring some of original theories in
\session{ZF-Constructible}, so these modifications are included as
part of our project.

It is to be noted that the proofs of the Forcing Theorems do not
require any extra replacement on the ground model; actually, they only need the 7
instances appearing in \isa{instances1{\uscore}fms} and
\isa{instances{\uscore}ground{\uscore}fms}.  But this seems not be
the case for Separation, at least by inspecting our formalization:
More instances holding in $M$ are needed 
as the complexity of $\phi$ grows. One point where this is apparent is
in the proof of Theorem~\ref{th:forcing-thms}(\ref{item:truth-lemma}),
that appears as the \isa{truth{\uscore}lemma} in our development; it
depends on \isa{truth{\uscore}lemma'} and
\isa{truth{\uscore}lemma{\uscore}Neg}, which explicitly invoke
\isa{separation{\uscore}ax}. In any case, our intended grounds
(v.g., the transitive collapse of countable elementary submodels of a
rank initial segment $V_\alpha$ or an $H(\kappa)$) all satisfy full
Separation.

\subsection{Models for $\CH$ and its negation}
\label{sec:models-ch-negation}

The statements of the existence of models of $\ZFC + \neg\CH$ and of
$\ZFC + \CH$  appear in our formalization as follows:

\begin{isabelle}
\isacommand{corollary}\isamarkupfalse%
\ ctm{\isacharunderscore}{\kern0pt}ZFC{\isacharunderscore}{\kern0pt}imp{\isacharunderscore}{\kern0pt}ctm{\isacharunderscore}{\kern0pt}not{\isacharunderscore}{\kern0pt}CH{\isacharcolon}{\kern0pt}\isanewline
\ \ \isakeyword{assumes}\isanewline
\ \ \ \ {\isachardoublequoteopen}M\ {\isasymapprox}\ {\isasymomega}{\isachardoublequoteclose}\ {\isachardoublequoteopen}Transset{\isacharparenleft}{\kern0pt}M{\isacharparenright}{\kern0pt}{\isachardoublequoteclose}\ {\isachardoublequoteopen}M\ {\isasymTurnstile}\ ZFC{\isachardoublequoteclose}\isanewline
\ \ \isakeyword{shows}\isanewline
\ \ \ \ {\isachardoublequoteopen}{\isasymexists}N{\isachardot}{\kern0pt}\isanewline
\ \ \ \ \ \ M\ {\isasymsubseteq}\ N\ {\isasymand}\ N\ {\isasymapprox}\ {\isasymomega}\ {\isasymand}\ Transset{\isacharparenleft}{\kern0pt}N{\isacharparenright}{\kern0pt}\ {\isasymand}\ N\ {\isasymTurnstile}\ ZFC\ {\isasymunion}\ {\isacharbraceleft}{\kern0pt}{\isasymcdot}{\isasymnot}{\isasymcdot}CH{\isasymcdot}{\isasymcdot}{\isacharbraceright}{\kern0pt}\ {\isasymand}\isanewline
\ \ \ \ \ \ {\isacharparenleft}{\kern0pt}{\isasymforall}{\isasymalpha}{\isachardot}{\kern0pt}\ Ord{\isacharparenleft}{\kern0pt}{\isasymalpha}{\isacharparenright}{\kern0pt}\ {\isasymlongrightarrow}\ {\isacharparenleft}{\kern0pt}{\isasymalpha}\ {\isasymin}\ M\ {\isasymlongleftrightarrow}\ {\isasymalpha}\ {\isasymin}\ N{\isacharparenright}{\kern0pt}{\isacharparenright}{\kern0pt}{\isachardoublequoteclose}
\end{isabelle}

\begin{isabelle}
\isacommand{corollary}\isamarkupfalse%
\ ctm{\isacharunderscore}{\kern0pt}ZFC{\isacharunderscore}{\kern0pt}imp{\isacharunderscore}{\kern0pt}ctm{\isacharunderscore}{\kern0pt}CH{\isacharcolon}{\kern0pt}\isanewline
\ \ \isakeyword{assumes}\isanewline
\ \ \ \ {\isachardoublequoteopen}M\ {\isasymapprox}\ {\isasymomega}{\isachardoublequoteclose}\ {\isachardoublequoteopen}Transset{\isacharparenleft}{\kern0pt}M{\isacharparenright}{\kern0pt}{\isachardoublequoteclose}\ {\isachardoublequoteopen}M\ {\isasymTurnstile}\ ZFC{\isachardoublequoteclose}\isanewline
\ \ \isakeyword{shows}\isanewline
\ \ \ \ {\isachardoublequoteopen}{\isasymexists}N{\isachardot}{\kern0pt}\isanewline
\ \ \ \ \ \ M\ {\isasymsubseteq}\ N\ {\isasymand}\ N\ {\isasymapprox}\ {\isasymomega}\ {\isasymand}\ Transset{\isacharparenleft}{\kern0pt}N{\isacharparenright}{\kern0pt}\ {\isasymand}\ N\ {\isasymTurnstile}\ ZFC\ {\isasymunion}\ {\isacharbraceleft}{\kern0pt}{\isasymcdot}CH{\isasymcdot}{\isacharbraceright}{\kern0pt}\ {\isasymand}\isanewline
\ \ \ \ \ \ {\isacharparenleft}{\kern0pt}{\isasymforall}{\isasymalpha}{\isachardot}{\kern0pt}\ Ord{\isacharparenleft}{\kern0pt}{\isasymalpha}{\isacharparenright}{\kern0pt}\ {\isasymlongrightarrow}\ {\isacharparenleft}{\kern0pt}{\isasymalpha}\ {\isasymin}\ M\ {\isasymlongleftrightarrow}\ {\isasymalpha}\ {\isasymin}\ N{\isacharparenright}{\kern0pt}{\isacharparenright}{\kern0pt}{\isachardoublequoteclose}
\end{isabelle}
where $\approx$ is equipotency, and the predicate \isa{Transset}
holds for
transitive sets. Both results are proved without using Choice.

As the excerpts indicate, these results are obtained as corollaries of
two theorems in which only a subset of the aforementioned
replacement instances are assumed of the ground model. We begin the
discussion of these stronger results by
considering extensions of ctms of fragments of $\ZF$.
\begin{isabelle}
\isacommand{theorem}\isamarkupfalse%
\ extensions{\isacharunderscore}{\kern0pt}of{\isacharunderscore}{\kern0pt}ctms{\isacharcolon}{\kern0pt}\isanewline
\ \ \isakeyword{assumes}\isanewline
\ \ \ \ {\isachardoublequoteopen}M\ {\isasymapprox}\ {\isasymomega}{\isachardoublequoteclose}\ {\isachardoublequoteopen}Transset{\isacharparenleft}{\kern0pt}M{\isacharparenright}{\kern0pt}{\isachardoublequoteclose}\isanewline
\ \ \ \ {\isachardoublequoteopen}M\ {\isasymTurnstile}\ {\isasymcdot}Z{\isasymcdot}\ {\isasymunion}\ {\isacharbraceleft}{\kern0pt}{\isasymcdot}Replacement{\isacharparenleft}{\kern0pt}p{\isacharparenright}{\kern0pt}{\isasymcdot}\ {\isachardot}{\kern0pt}\ p\ {\isasymin}\ overhead{\isacharbraceright}{\kern0pt}{\isachardoublequoteclose}\isanewline
\ \ \ \ {\isachardoublequoteopen}{\isasymPhi}\ {\isasymsubseteq}\ formula{\isachardoublequoteclose}\isanewline%
\ \ \ \ {\isachardoublequoteopen}M\ {\isasymTurnstile}\ {\isacharbraceleft}{\kern0pt}\ {\isasymcdot}Replacement{\isacharparenleft}{\kern0pt}ground{\isacharunderscore}{\kern0pt}repl{\isacharunderscore}{\kern0pt}fm{\isacharparenleft}{\kern0pt}{\isasymphi}{\isacharparenright}{\kern0pt}{\isacharparenright}{\kern0pt}{\isasymcdot}\ {\isachardot}{\kern0pt}\ {\isasymphi}\ {\isasymin}\ {\isasymPhi}{\isacharbraceright}{\kern0pt}{\isachardoublequoteclose}\isanewline
\ \ \isakeyword{shows}\isanewline
\ \ \ \ {\isachardoublequoteopen}{\isasymexists}N{\isachardot}{\kern0pt}\isanewline
\ \ \ \ \ \ M\ {\isasymsubseteq}\ N\ {\isasymand}\ N\ {\isasymapprox}\ {\isasymomega}\ {\isasymand}\ Transset{\isacharparenleft}{\kern0pt}N{\isacharparenright}{\kern0pt}\ {\isasymand}\ M{\isasymnoteq}N\ {\isasymand}\isanewline
\ \ \ \ \ \ {\isacharparenleft}{\kern0pt}{\isasymforall}{\isasymalpha}{\isachardot}{\kern0pt}\ Ord{\isacharparenleft}{\kern0pt}{\isasymalpha}{\isacharparenright}{\kern0pt}\ {\isasymlongrightarrow}\ {\isacharparenleft}{\kern0pt}{\isasymalpha}\ {\isasymin}\ M\ {\isasymlongleftrightarrow}\ {\isasymalpha}\ {\isasymin}\ N{\isacharparenright}{\kern0pt}{\isacharparenright}{\kern0pt}\ {\isasymand}\isanewline
\ \ \ \ \ \ {\isacharparenleft}{\kern0pt}{\isacharparenleft}{\kern0pt}M{\isacharcomma}{\kern0pt}\ {\isacharbrackleft}{\kern0pt}{\isacharbrackright}{\kern0pt}{\isasymTurnstile}\ {\isasymcdot}AC{\isasymcdot}{\isacharparenright}{\kern0pt}\ {\isasymlongrightarrow}\ N{\isacharcomma}{\kern0pt}\ {\isacharbrackleft}{\kern0pt}{\isacharbrackright}{\kern0pt}\ {\isasymTurnstile}\ {\isasymcdot}AC{\isasymcdot}{\isacharparenright}{\kern0pt}\ {\isasymand}\isanewline
\ \ \ \ \ \ N\ {\isasymTurnstile}\ {\isasymcdot}Z{\isasymcdot}\ {\isasymunion}\ {\isacharbraceleft}{\kern0pt}\ {\isasymcdot}Replacement{\isacharparenleft}{\kern0pt}{\isasymphi}{\isacharparenright}{\kern0pt}{\isasymcdot}\ {\isachardot}{\kern0pt}\ {\isasymphi}\ {\isasymin}\ {\isasymPhi}{\isacharbraceright}{\kern0pt}{\isachardoublequoteclose}
\end{isabelle}

Here, the 7-element set \isa{overhead} is enough to construct a proper
extension. It is  the union of
\isa{instances{\isadigit{1}}{\isacharunderscore}{\kern0pt}fms} and
\isa{instances{\isacharunderscore}{\kern0pt}ground{\isacharunderscore}{\kern0pt}fms}.
Also,
\isa{{\isasymcdot}Z{\isasymcdot}} denotes Zermelo set theory and one
can use the parameter $\Phi$ to ensure those replacement instances in the extension.

In the
next theorem, the relevant set of formulas is
\isa{overhead{\isacharunderscore}{\kern0pt}notCH}, defined above in
Section~\ref{sec:repl-instances}, and \isa{ZC} denotes Zermelo set
theory plus Choice:

\begin{isabelle}
\isacommand{theorem}\isamarkupfalse%
\ ctm{\isacharunderscore}{\kern0pt}of{\isacharunderscore}{\kern0pt}not{\isacharunderscore}{\kern0pt}CH{\isacharcolon}{\kern0pt}\isanewline
\ \ \isakeyword{assumes}\isanewline
\ \ \ \ {\isachardoublequoteopen}M\ {\isasymapprox}\ {\isasymomega}{\isachardoublequoteclose}\ {\isachardoublequoteopen}Transset{\isacharparenleft}{\kern0pt}M{\isacharparenright}{\kern0pt}{\isachardoublequoteclose}\isanewline
\ \ \ \ {\isachardoublequoteopen}M\ {\isasymTurnstile}\ ZC\ {\isasymunion}\ {\isacharbraceleft}{\kern0pt}{\isasymcdot}Replacement{\isacharparenleft}{\kern0pt}p{\isacharparenright}{\kern0pt}{\isasymcdot}\ {\isachardot}{\kern0pt}\ p\ {\isasymin}\ overhead{\isacharunderscore}{\kern0pt}notCH{\isacharbraceright}{\kern0pt}{\isachardoublequoteclose}\isanewline
\ \ \ \ {\isachardoublequoteopen}{\isasymPhi}\ {\isasymsubseteq}\ formula{\isachardoublequoteclose}\isanewline
\ \ \ \ {\isachardoublequoteopen}M\ {\isasymTurnstile}\ {\isacharbraceleft}{\kern0pt}\ {\isasymcdot}Replacement{\isacharparenleft}{\kern0pt}ground{\isacharunderscore}{\kern0pt}repl{\isacharunderscore}{\kern0pt}fm{\isacharparenleft}{\kern0pt}{\isasymphi}{\isacharparenright}{\kern0pt}{\isacharparenright}{\kern0pt}{\isasymcdot}\ {\isachardot}{\kern0pt}\ {\isasymphi}\ {\isasymin}\ {\isasymPhi}{\isacharbraceright}{\kern0pt}{\isachardoublequoteclose}\isanewline
\ \ \isakeyword{shows}\isanewline
\ \ \ \ {\isachardoublequoteopen}{\isasymexists}N{\isachardot}{\kern0pt}\isanewline
\ \ \ \ \ \ M\ {\isasymsubseteq}\ N\ {\isasymand}\ N\ {\isasymapprox}\ {\isasymomega}\ {\isasymand}\ Transset{\isacharparenleft}{\kern0pt}N{\isacharparenright}{\kern0pt}\ {\isasymand}\isanewline
\ \ \ \ \ \ N\ {\isasymTurnstile}\ ZC\ {\isasymunion}\ {\isacharbraceleft}{\kern0pt}{\isasymcdot}{\isasymnot}{\isasymcdot}CH{\isasymcdot}{\isasymcdot}{\isacharbraceright}{\kern0pt}\ {\isasymunion}\ {\isacharbraceleft}{\kern0pt}\ {\isasymcdot}Replacement{\isacharparenleft}{\kern0pt}{\isasymphi}{\isacharparenright}{\kern0pt}{\isasymcdot}\ {\isachardot}{\kern0pt}\ {\isasymphi}\ {\isasymin}\ {\isasymPhi}{\isacharbraceright}{\kern0pt}\ {\isasymand}\isanewline
\ \ \ \ \ \ {\isacharparenleft}{\kern0pt}{\isasymforall}{\isasymalpha}{\isachardot}{\kern0pt}\ Ord{\isacharparenleft}{\kern0pt}{\isasymalpha}{\isacharparenright}{\kern0pt}\ {\isasymlongrightarrow}\ {\isacharparenleft}{\kern0pt}{\isasymalpha}\ {\isasymin}\ M\ {\isasymlongleftrightarrow}\ {\isasymalpha}\ {\isasymin}\ N{\isacharparenright}{\kern0pt}{\isacharparenright}{\kern0pt}{\isachardoublequoteclose}
\end{isabelle}

Finally, \isa{overhead{\isacharunderscore}{\kern0pt}CH} is the union
of \isa{overhead{\isacharunderscore}{\kern0pt}notCH} with the $\DC$
instance \isa{dc{\uscore}abs{\uscore}fm}:
\begin{isabelle}
\isacommand{theorem}\isamarkupfalse%
\ ctm{\isacharunderscore}{\kern0pt}of{\isacharunderscore}{\kern0pt}CH{\isacharcolon}{\kern0pt}\isanewline
\ \ \isakeyword{assumes}\isanewline
\ \ \ \ {\isachardoublequoteopen}M\ {\isasymapprox}\ {\isasymomega}{\isachardoublequoteclose}\ {\isachardoublequoteopen}Transset{\isacharparenleft}{\kern0pt}M{\isacharparenright}{\kern0pt}{\isachardoublequoteclose}\isanewline
\ \ \ \ {\isachardoublequoteopen}M\ {\isasymTurnstile}\ ZC\ {\isasymunion}\ {\isacharbraceleft}{\kern0pt}{\isasymcdot}Replacement{\isacharparenleft}{\kern0pt}p{\isacharparenright}{\kern0pt}{\isasymcdot}\ {\isachardot}{\kern0pt}\ p\ {\isasymin}\ overhead{\isacharunderscore}{\kern0pt}CH{\isacharbraceright}{\kern0pt}{\isachardoublequoteclose}\isanewline
\ \ \ \ {\isachardoublequoteopen}{\isasymPhi}\ {\isasymsubseteq}\ formula{\isachardoublequoteclose}\isanewline
\ \ \ \ {\isachardoublequoteopen}M\ {\isasymTurnstile}\ {\isacharbraceleft}{\kern0pt}\ {\isasymcdot}Replacement{\isacharparenleft}{\kern0pt}ground{\isacharunderscore}{\kern0pt}repl{\isacharunderscore}{\kern0pt}fm{\isacharparenleft}{\kern0pt}{\isasymphi}{\isacharparenright}{\kern0pt}{\isacharparenright}{\kern0pt}{\isasymcdot}\ {\isachardot}{\kern0pt}\ {\isasymphi}\ {\isasymin}\ {\isasymPhi}{\isacharbraceright}{\kern0pt}{\isachardoublequoteclose}\isanewline
\ \ \isakeyword{shows}\isanewline
\ \ \ \ {\isachardoublequoteopen}{\isasymexists}N{\isachardot}{\kern0pt}\isanewline
\ \ \ \ \ \ M\ {\isasymsubseteq}\ N\ {\isasymand}\ N\ {\isasymapprox}\ {\isasymomega}\ {\isasymand}\ Transset{\isacharparenleft}{\kern0pt}N{\isacharparenright}{\kern0pt}\ {\isasymand}\isanewline
\ \ \ \ \ \ N\ {\isasymTurnstile}\ ZC\ {\isasymunion}\ {\isacharbraceleft}{\kern0pt}{\isasymcdot}CH{\isasymcdot}{\isacharbraceright}{\kern0pt}\ {\isasymunion}\ {\isacharbraceleft}{\kern0pt}\ {\isasymcdot}Replacement{\isacharparenleft}{\kern0pt}{\isasymphi}{\isacharparenright}{\kern0pt}{\isasymcdot}\ {\isachardot}{\kern0pt}\ {\isasymphi}\ {\isasymin}\ {\isasymPhi}{\isacharbraceright}{\kern0pt}\ {\isasymand}\isanewline
\ \ \ \ \ \ {\isacharparenleft}{\kern0pt}{\isasymforall}{\isasymalpha}{\isachardot}{\kern0pt}\ Ord{\isacharparenleft}{\kern0pt}{\isasymalpha}{\isacharparenright}{\kern0pt}\ {\isasymlongrightarrow}\ {\isacharparenleft}{\kern0pt}{\isasymalpha}\ {\isasymin}\ M\ {\isasymlongleftrightarrow}\ {\isasymalpha}\ {\isasymin}\ N{\isacharparenright}{\kern0pt}{\isacharparenright}{\kern0pt}{\isachardoublequoteclose}
\end{isabelle}

\section{Related work}
\label{sec:related-work}

There is another formalization of forcing in Lean by Han and van
Doorn, under the name \emph{Flypitch} \cite{han_et_al:LIPIcs:2019:11074,DBLP:conf/cpp/HanD20}. When our project started, we
were unaware of this initiative, and the same as them, we were deeply
influenced by Wiedijk's list of 100 theorems \cite{Formalizing100}.

Many aspects make their formalization different from ours. Their
presentation of the mathematics is somewhat more elegant and cohesive,
since they go for the Boolean valued approach; they also  set up the
calculus of first-order logic, and \emph{en route} to forcing they
formalized the basic model theory of Boolean valued models and Gödel's
Completeness Theorem. They also provided the treatment of the regular
open algebra, and the general version of the Delta System
Lemma. Putting this together they readily obtain a proof that
$\ZFC\nvdash \CH$ \cite{han_et_al:LIPIcs:2019:11074} and after
formalizing collapse forcing they show  $\ZFC\nvdash \neg\CH$
\cite[Sect.~5.6]{DBLP:conf/cpp/HanD20}.

It should be emphasized, however, that the Flypitch project was
carried out assuming a rather strong metatheory.
Carneiro \cite{carneiro-ms-thesis} reports that Werner's results in
\cite{10.5555/645869.668660} can be adapted to show that the base
logic  of Lean (restricted to $n$
type universes) proves the consistency of $\ZFC$ plus $n$
inaccessibles. Han and van Doorn did use universes in their
implementation; for instance, ordinals are “defined as equivalence
classes of (well-ordered) types, [\dots] one
universe level higher than the types used to construct them”
\cite{han_et_al:LIPIcs:2019:11074}. It is not clear to us if they are
able to avoid such strength: At least, $\Con(\ZFC)$ is provable in
their context. %
On a lesser note, in order to prove $\AC$ in the generic extension,
Flypitch requires choice in the metatheory
\cite[p.~11]{han_et_al:LIPIcs:2019:11074}, while our formalization works
entirely in $\ZF$. %

This is perhaps an appropriate time to insist that we have \emph{not}
formalized the relative consistency of $\neg\CH$, and we are actually
not aiming for that in the short term. In our context, going for the
plain consistency result (that is, proving the arithmetical statement
$\Con(ZFC) \implies \Con(ZFC + \neg\CH)$ over a weak base theory)
seems off the mark, since our base is fixed (and essentially equivalent to $\ZF$). %
If we intended to do so, we would need a development of a weak base with a
new type of (external) formulas having the required
induction principles, plus an interface from this new type to
Isabelle/ZF. This might allow us to take the standard route to use ctms for a
relative consistency proof (through a version of the Reflection
Theorem proved by induction on the external formulas).

Back to the comparison, we believe that our formalization using the ctm approach
over Isabelle/ZF might be more appealing to set-theorists because of the
type-theoretic machinery used in Flypitch, and %
since absoluteness grants us extra naturality. This last point may also
be illustrated by the treatment of ordinals; in our formalization, as
it is expected intuitively, the following are equivalent for every $x$
in a ctm $M$:
\begin{itemize}
\item $\isa{Ord}(x)$;
\item $\isa{ordinal}(\isa{\#\#}M,x)$ (the relational
  definition relativized to $M$);
\item
  $M,[x] \models {\cdot}0\isa{ is ordinal}{\cdot}$.
\end{itemize}
where ${\cdot}n\isa{ is ordinal}{\cdot}$ is the code for the
appropriate first-order formula  ($0$ is a de
Bruijn index above!). In contrast, Han and van Doorn require an
injection from the ordinals of the corresponding type universe into
their encoding of a model of $\ZF$, and a further necessary injection into the Boolean
valued model using checks---this last step obviously appears in our
presentation, but the $\val$ function used to construct $M[G]$ will turn
check-names into the corresponding argument, as expected.

This faithfulness to set-theoretical practice does not come for
free. Recursive constructions and inductive definitions are far easier
to perform in the Calculus of Inductive Constructions on which Lean
is based, and in Isabelle/ZF are rather cumbersome. Also, a
typed discipline provides aid to write succinctly and many assumptions
are satisfied by mere notation. To be clear, those benefits come from
doing set theory in a non set-theoretical language. On the other
hand, Isar proofs, as the one shown in
Sec.~\ref{sec:sample-formal-proof}, are easier to understand
than the language of tactics of Lean.

A sweet spot combining the best of both worlds is to be found on
developments in Isabelle/HOL based on the AFP entry
\session{ZFC\_in\_HOL} by Paulson \cite{ZFC_in_HOL-AFP}. There is a
range of results in combinatorics and other set-theoretical material
that was swiftly formalized in this setting: Erdős-Milner
partition theorem \cite{2021arXiv210411613P}, Nash-Williams theorem and
Larson's $\forall k\ \omega^{\omega}\longrightarrow(\omega^\omega,k)$
\cite{doi:10.1080/10586458.2021.1980464}, Design Theory
\cite{10.1007/978-3-030-81097-9_1}, and Wetzel's problem
\cite{2022arXiv220503159P}; this last paper describes in a brief and
clear way the convenience of the interaction with Isabelle/HOL
(which is also paid in consistency-strength currency
\cite[Sect.~3]{DBLP:conf/ictac/Obua06}).

Concerning the minimum amount of Replacement needed to construct
forcing extensions, only recently we learned about Mathias' work on the
subject (for which a summary is offered in
\cite[Sect.~6]{kanamori-mathias}). In
\cite[Sect.~1]{mathias:hal-01188043}, models $M$ of Zermelo set theory
are constructed for which each of the inclusions $M\subseteq M[G]$ and
$M\supseteq M[G]$ fail, where the poset $\pP$ is the trivial $\{\1
\}$. In one of them, $K$, we have $\om\in K \sm K[G]$, hence the
ordinals do not coincide.

Also, in the reference \cite{mathias-provident}, a reasonably minimal
fragment \textsc{Prov} of $\ZF$ that allows to do set forcing is
identified, and transitive sets satisfying it are called
\emph{provident}. Existence of rank and of transitive closure are
implied by \textsc{Prov}; hence their appearance in our list seems
justified. Nevertheless, \textsc{Prov} is far weaker than the
fragments of $\ZF$ considered here, since it is a restriction of
Kripke-Platek set theory, and thus it does not include neither
Powerset nor full Separation. The detailed theory of provident sets is
developed in \cite{mathias-bowler-gentle}. %

\section{Some lessons}\label{sec:lessons}

We want to finish this report by gathering some of the conclusions we
reached after the experience of formalizing the basics of forcing in a
proof assistant.

\subsection{Aims of a formalization and planning}
\label{sec:aims-formalization-planning}

We believe that in every project of formalization of mathematics,
there is a tension between the haste to verify the target results and
the need to obtain a readable, albeit extremely detailed, corpus of
statements and proofs. This tension is mirrored in two different
purposes of formalization: Developing new mathematics from scratch and
producing verified results on the way, versus verifying and
documenting material that has already been produced on paper.

Our present project clearly belongs to the second category, so we
prioritized trying to obtain formal proofs that mimicked standard
prose (as can be seen in the sample proof in
Section~\ref{sec:sample-formal-proof}). We feel that the Isar language
provided with Isabelle has the right balance between elegance and
efficacy. Another crucial aspect to achieve this goal is the level of
detail of a blueprint for the formalization.

We must however confess
that we learned many of the subtleties of Isabelle in the making, and
many engineering decisions were also taken before it was clear the
precise way things would develop in the future.
From time to time, we were frenziedly trying to
get the results formalized, going beyond what we had planned.

As a result from this, some design choices that seemed reasonable at
first were proved to be inconvenient. For instance, we should had
better used predicates (of type $\tyi\fun\tyi\fun\tyo$) for the
forcing posets' order relations; this is the way they
are presented in the \session{Delta\_System\_Lemma} session. A similar
problem, which can be traced to our reading of Kunen's suggestion on
how to formalize the forcing relation \cite[p.~260]{kunen2011set},
is that we require the forcing poset to be an element of $M$,
so the present infrastructure does not allow class forcing out of the
box. (The latter change seems to be rather straightforward, but the
former does not.)

Upon reaching the final stage of the project, we decided to go for the minimal
set of definitions and versions of lemmas that were needed to obtain
our target results. For example,
\begin{itemize}
\item
  we only proved the Delta System Lemma for $\aleph_1$-sized families;
  thus limiting us to the case of the $\aleph_1$-chain condition, and avoiding the relativization
  of the material on cofinalities \cite{Delta_System_Lemma-AFP};
\item we showed preservation of sequences by considering countably
  closed forcings (in fact, we formalized the bare minimum requirement
  of being $({<}\omega{+}1)$-closed, that is, closure under
  $\omega$-sequences and not $\delta$-sequences for every countable $\delta$).
\end{itemize}
In doing this we went against the conventional wisdom that one should formalize the most general version of the
results available. Another shortcut we took was to
simplify some proofs by appealing to the countability of the ground
model; this is the case of
\isa{definition{\uscore}of{\uscore}forces} and the result on forcing
values of a function.

\subsection{How to believe in the formalization}
\label{sec:how-believe-formalization}

This is a rather tricky question, that was addressed by Pollack in his
\cite{MR1686867}. There is little point to discuss that, after an
assistant has accepted some input successfully, \emph{some
mathematics} has been formally verified. What might not be apparent is
if the claimed theorems are indeed the results that have been
checked. One key aspect of this is the logical foundation of the
assistant (Section~\ref{sec:isabelle-metalogic-meta}). But the weakest
link in the chain is the laying down of definitions building up to the
concepts needed to state the target results.%
\footnote{%
Another related aspect,
concerning the way results are printed and parsed by assistant versus
their internal meaning, was studied by Wiedijk \cite{zbMATH06319597}.}

We took care of this matter by providing, as an entry point for our
whole development, the theory \theory{Definitions\_Main} in which a
path from some the fundamental concepts from Isabelle/ZF reaching to our main theorems
is expounded. Cross-references to major milestones (which can be
navigated by using Isabelle) are provided there. A curated version can
be found as Appendix~\ref{sec:definitions_main} to this paper.

Frequently, we formalized material by directly typing the proof we
knew by heart, and in so doing we assumed that some definitions
accommodated some of our preconceptions.
It is significant that in a few such occasions, we were doubly
surprised by the fact that some supposedly trivial lemma would not go
through, because the definitions addressed something different (think
of \emph{restriction} of a function to a set versus that of a
relation), and also that we were able to prove the adjacent
results. The takeaway is that intuition may drive proofs
even if you are not working on what you think you are.

A final aspect on this topic concerns automated methods. In the
Introduction we hinted at the fact that a proof can actually be
\emph{obscured} by automation. Specifically, proof steps that were solved automatically give no
information for someone who wants to understand the details of the
argument; by the same token, automatic methods might silently exploit
inconsistencies in the definitions, and this will only be apparent in
a later stage of the development.

\subsection{Bureaucracy and scale factors}
\label{sec:bureaucracy-scale-factors}

It is noteworthy that although the “math” of the construction of a
model of $\ZFC+\neg\CH$ was already in place by the end of November
2020, it was only 9 months later that we were able to finish the
formalization of that result. The missing pieces were essentially
bureaucracy. Some of the material filed under this category comprises:
\begin{itemize}
\item permutation of indices and calculation of arities of
  internalized formulas;
\item proving that certain constructions belong to the relevant
  models;
\item (required for the above) showing that particular instances of
  separation and replacement hold in the ground model.
\end{itemize}

Some of those proofs were almost copy-pasted once and again with minor
variants; this would usually be relegated to some function in the
meta-language, but we were unable to do this due to our limitations in
programming Isabelle/ML\footnote{%
  On the other hand, our inability to automate proofs of replacement
  instances paved the way for identifying which were the ones needed for
  forcing!}.

Nevertheless, experience  in software engineering is invaluable in
large projects like the present one. For example, it is
(mathematically) misleading when automatic tools (\isa{simp},
\isa{auto}, etc) stop working just because of the sheer size of the
goal (v.g., the same statement with 7 variables succeeds but with 8
variables does not). Scale issues are very easily disregarded in the
abstract but, as a colorful example, the formula $\forceisa({\cdot}0\in 1{\cdot})$
can be explicitly printed by Isabelle2021-1 (it spans nearly 20k symbols), but
$\forceisa({\cdot}\neg{\cdot}\neg{\cdot}0\in 1{\cdot}{\cdot}{\cdot})$
can not.
(For this particular example, one reason for the blow-up in
length comes from the encoding of negation using the \isatt{Nand}
connective, but otherwise the formulas grow bigger by the sheer
complexity of the definition of the forcing relation.)

Another point where computer science expertise was a prime asset was
the very definition of $\forceisa$. As a
proof of concept, one of us tried to obtain its definition by
using formula synthesis exclusively, which was supposed to be as
trivial as in the usual
mathematical development (similarly to the case of
Equation~(\ref{eq:sats_big_union_fm})).
But in fact, some early minor mistake
rendered the whole effort useless. We then turned to a more
informed programming discipline, which involved decomposing the
definition in stages, each of which was checked for correctness, and
in that way we were able to reach our objective.

\section{Future directions}
\label{sec:conclusion}

There are many possibilities for further work starting from this
formalization. We will mention just a few.

Obvious missing pieces would be proving the standard properties of
general Cohen posets $\Fn_\kappa(I,J)$, and to modify the core
definitions to allow for class forcing. We would also like to try the
Boolean valued approach to compare the (un)ease of formalization using
Isabelle/ZF.

Another desirable goal is to construct transitive set models of $\ZFC$
from a large cardinal. There is some work to be done for that: Even the
definition of inaccessibles, and of the transitive collapse (for that
matter), are still missing.

As we did for Replacement, we would like to pinpoint an (almost)
minimal set of instances of
Separation needed to use forcing. A necessary ingredient will
certainly be an implementation of Gödel operations
\cite[Thm.~13.4]{Jech_Millennium}.

In our previous landmark \cite{2020arXiv200109715G}, we contributed
with some modifications to \session{ZF-Constructible}; this is now
part of the official Isabelle distribution.  We intend ask Isabelle
maintainers to consider the more modular versions of some of those
theories that we are presenting in this project.

As final words about our journey, we believe that, as in mathematics
in general, the experience of working in a formal environment can be
daunting, but at the same time extremely rewarding: The feeling of
accomplishment after seeing your own writings validated beyond doubt is
in some ways comparable to that of finding a proof of an important
lemma. It also allows subtly different ways of reasoning (with their
own merits and pitfalls---it is easy to forget how easy a proof on
paper is once you are fully engaged in directing your assistant). We
hope that at some point these experiences are shared by our community
at large.

\section*{Acknowledgments}
\label{sec:acknowledgments}
We are indebted to Timothy Chow for insightful comments on a draft of
this paper. We also thank Z.~Vidnyánszky for some discussions where the idea of
building the summary in \theory{Definitions\_Main} took shape.
Johan Commelin provided useful comments on LTE.
We are very grateful
for the anonymous referee's comments and pointers, which enhanced our presentation.
Mikhail Mandrykin helped us with technical assistance regarding the
rendering of formulas in the Isabelle/jEdit interface, which improved
our presentation of this material during the Latin American Congress
of Mathematicians (CLAM 2021); we extend our gratitude to the rest of
the Isabelle community, for their kind support through the Users mailing
list. We also want to warmly thank Larry Paulson for his work and
encouragement. Last, and definitely not least, we are deeply grateful
to Ken Kunen for his
inspiring expositions. This formalization is a tiny homage to his
memory.

\providecommand{\noopsort}[1]{}
\begin{small}\end{small}

\appendix

\section{Main definitions of the formalization}\label{sec:definitions_main}

This section, which appears almost verbatim as
the theory \theory{Definitions\_Main} in \cite{Independence_CH-AFP},
might be considered as the bare minimum reading requisite to
trust that our development indeed formalizes the theory of
forcing.

The reader trusting
all the libraries on which our development is based, might jump
directly to Section~\ref{sec:relative-arith}, which treats relative
cardinal arithmetic as implemented in
\isa{T{\kern0pt}r{\kern0pt}a{\kern0pt}n{\kern0pt}s{\kern0pt}i{\kern0pt}t{\kern0pt}i{\kern0pt}v{\kern0pt}e{\kern0pt}{\char`\_}{\kern0pt}M{\kern0pt}o{\kern0pt}d{\kern0pt}e{\kern0pt}l{\kern0pt}s{\kern0pt}}. But in case one wants to dive deeper, the
following sections treat some basic concepts of the ZF logic
(Section~\ref{sec:def-main-ZF}) and in the
\session{ZF-Constructible} library (Section~\ref{sec:def-main-relative})
on which our definitions are built.

\subsection{ZF\label{sec:def-main-ZF}%
}
For the basic logic ZF we restrict ourselves to just a few
concepts (for its axioms, consult Appendix~\ref{appendix:axioms}).
\begin{isabelle}%
bij{\isacharparenleft}{\kern0pt}A{\isacharcomma}{\kern0pt}\ B{\isacharparenright}{\kern0pt}\ {\isasymequiv}\isanewline
{\isacharbraceleft}{\kern0pt}f\ {\isasymin}\ A\ {\isasymrightarrow}\ B\ {\isachardot}{\kern0pt}\ {\isasymforall}w{\isasymin}A{\isachardot}{\kern0pt}\ {\isasymforall}x{\isasymin}A{\isachardot}{\kern0pt}\ f\ {\isacharbackquote}{\kern0pt}\ w\ {\isacharequal}{\kern0pt}\ f\ {\isacharbackquote}{\kern0pt}\ x\ {\isasymlongrightarrow}\ w\ {\isacharequal}{\kern0pt}\ x{\isacharbraceright}{\kern0pt}\ {\isasyminter}\isanewline
{\isacharbraceleft}{\kern0pt}f\ {\isasymin}\ A\ {\isasymrightarrow}\ B\ {\isachardot}{\kern0pt}\ {\isasymforall}y{\isasymin}B{\isachardot}{\kern0pt}\ {\isasymexists}x{\isasymin}A{\isachardot}{\kern0pt}\ f\ {\isacharbackquote}{\kern0pt}\ x\ {\isacharequal}{\kern0pt}\ y{\isacharbraceright}{\kern0pt}%
\end{isabelle}%
\begin{isabelle}%
A\ {\isasymapprox}\ B\ {\isasymequiv}\ {\isasymexists}f{\isachardot}{\kern0pt}\ f\ {\isasymin}\ bij{\isacharparenleft}{\kern0pt}A{\isacharcomma}{\kern0pt}\ B{\isacharparenright}{\kern0pt}%
\end{isabelle}%
\begin{isabelle}%
Transset{\isacharparenleft}{\kern0pt}i{\isacharparenright}{\kern0pt}\ {\isasymequiv}\ {\isasymforall}x{\isasymin}i{\isachardot}{\kern0pt}\ x\ {\isasymsubseteq}\ i%
\end{isabelle}%
\begin{isabelle}%
Ord{\isacharparenleft}{\kern0pt}i{\isacharparenright}{\kern0pt}\ {\isasymequiv}\ Transset{\isacharparenleft}{\kern0pt}i{\isacharparenright}{\kern0pt}\ {\isasymand}\ {\isacharparenleft}{\kern0pt}{\isasymforall}x{\isasymin}i{\isachardot}{\kern0pt}\ Transset{\isacharparenleft}{\kern0pt}x{\isacharparenright}{\kern0pt}{\isacharparenright}{\kern0pt}%
\end{isabelle}%
\begin{isabelle}%
i\ {\isacharless}{\kern0pt}\ j\ {\isasymequiv}\ i\ {\isasymin}\ j\ {\isasymand}\ Ord{\isacharparenleft}{\kern0pt}j{\isacharparenright}{\kern0pt}\isasep\isanewline%
i\ {\isasymle}\ j\ {\isasymlongleftrightarrow}\ i\ {\isacharless}{\kern0pt}\ j\ {\isasymor}\ {\isacharparenleft}i\ {\isacharequal}{\kern0pt}\ j\ {\isasymand}\ Ord{\isacharparenleft}{\kern0pt}j{\isacharparenright}{\isacharparenright}{\kern0pt}%
\end{isabelle}%
With the concepts of empty set and successor in place,%

\begin{isabelle}
\isacommand{lemma}
\ empty{\uscore}{\kern0pt}def{\isacharprime}{\kern0pt}{\isacharcolon}{\kern0pt}\ {\isachardoublequoteopen}{\isasymforall}x{\isachardot}{\kern0pt}\ x\ {\isasymnotin}\ {\isadigit{0}}{\isachardoublequoteclose}%
\isanewline
\isacommand{lemma}
\ succ{\uscore}{\kern0pt}def{\isacharprime}{\kern0pt}{\isacharcolon}{\kern0pt}\ {\isachardoublequoteopen}succ{\isacharparenleft}{\kern0pt}i{\isacharparenright}{\kern0pt}\ {\isacharequal}{\kern0pt}\ i\ {\isasymunion}\ {\isacharbraceleft}{\kern0pt}i{\isacharbraceright}{\kern0pt}{\isachardoublequoteclose}%
\end{isabelle}
we can define the set of natural numbers \isa{{\isasymomega}}. In the
sources, it is  defined as a fixpoint, but here we just write
its characterization as the first limit ordinal.%
\begin{isabelle}%
Ord{\isacharparenleft}{\kern0pt}{\isasymomega}{\isacharparenright}{\kern0pt}\ {\isasymand}\ {\isadigit{0}}\ {\isacharless}{\kern0pt}\ {\isasymomega}\ {\isasymand}\ {\isacharparenleft}{\kern0pt}{\isasymforall}y{\isachardot}{\kern0pt}\ y\ {\isacharless}{\kern0pt}\ {\isasymomega}\ {\isasymlongrightarrow}\ succ{\isacharparenleft}{\kern0pt}y{\isacharparenright}{\kern0pt}\ {\isacharless}{\kern0pt}\ {\isasymomega}{\isacharparenright}{\kern0pt}\isasep\isanewline%
Ord{\isacharparenleft}{\kern0pt}i{\isacharparenright}{\kern0pt}\ {\isasymand}\ {\isadigit{0}}\ {\isacharless}{\kern0pt}\ i\ {\isasymand}\ {\isacharparenleft}{\kern0pt}{\isasymforall}y{\isachardot}{\kern0pt}\ y\ {\isacharless}{\kern0pt}\ i\ {\isasymlongrightarrow}\ succ{\isacharparenleft}{\kern0pt}y{\isacharparenright}{\kern0pt}\ {\isacharless}{\kern0pt}\ i{\isacharparenright}{\kern0pt}\ {\isasymLongrightarrow}\ {\isasymomega}\ {\isasymle}\ i%
\end{isabelle}%
Then, addition and predecessor on \isa{{\isasymomega}} are inductively
characterized as follows:%
\begin{isabelle}%
m\ {\isacharplus}{\kern0pt}\isactrlsub {\isasymomega}\ succ{\isacharparenleft}{\kern0pt}n{\isacharparenright}{\kern0pt}\ {\isacharequal}{\kern0pt}\ succ{\isacharparenleft}{\kern0pt}m\ {\isacharplus}{\kern0pt}\isactrlsub {\isasymomega}\ n{\isacharparenright}{\kern0pt}\isasep\isanewline%
m\ {\isasymin}\ {\isasymomega}\ {\isasymLongrightarrow}\ m\ {\isacharplus}{\kern0pt}\isactrlsub {\isasymomega}\ {\isadigit{0}}\ {\isacharequal}{\kern0pt}\ m\isasep\isanewline\isanewline%
pred{\isacharparenleft}{\kern0pt}{\isadigit{0}}{\isacharparenright}{\kern0pt}\ {\isacharequal}{\kern0pt}\ {\isadigit{0}}\isasep\isanewline%
pred{\isacharparenleft}{\kern0pt}succ{\isacharparenleft}{\kern0pt}y{\isacharparenright}{\kern0pt}{\isacharparenright}{\kern0pt}\ {\isacharequal}{\kern0pt}\ y%
\end{isabelle}%
Lists on a set \isa{A} can be characterized by being
recursively generated from the empty list \isa{{\isacharbrackleft}{\kern0pt}{\isacharbrackright}{\kern0pt}} and the
operation \isa{Cons} that adds a new element to the left end;
the induction theorem for them shows that the characterization is
“complete”. (Mind the
\isa{\isasymlbrakk P; Q\isasymrbrakk\ \isasymLongrightarrow\ R}
abbreviation for
\isa{P\ \isasymLongrightarrow\ Q\ \isasymLongrightarrow\ R}.)

\begin{isabelle}%
{\isacharbrackleft}{\kern0pt}{\isacharbrackright}{\kern0pt}\ {\isasymin}\ list{\isacharparenleft}{\kern0pt}A{\isacharparenright}{\kern0pt}\isasep\isanewline%
{\isasymlbrakk}a\ {\isasymin}\ A{\isacharsemicolon}{\kern0pt}\ l\ {\isasymin}\ list{\isacharparenleft}{\kern0pt}A{\isacharparenright}{\kern0pt}{\isasymrbrakk}\ {\isasymLongrightarrow}\ Cons{\isacharparenleft}{\kern0pt}a{\isacharcomma}{\kern0pt}\ l{\isacharparenright}{\kern0pt}\ {\isasymin}\ list{\isacharparenleft}{\kern0pt}A{\isacharparenright}{\kern0pt}\isasep\isanewline\isanewline%
{\isasymlbrakk}x\ {\isasymin}\ list{\isacharparenleft}{\kern0pt}A{\isacharparenright}{\kern0pt}{\isacharsemicolon}{\kern0pt}\ P{\isacharparenleft}{\kern0pt}{\isacharbrackleft}{\kern0pt}{\isacharbrackright}{\kern0pt}{\isacharparenright}{\kern0pt}{\isacharsemicolon}{\kern0pt}\ {\isasymAnd}a\ l{\isachardot}{\kern0pt}\ {\isasymlbrakk}a\ {\isasymin}\ A{\isacharsemicolon}{\kern0pt}\ l\ {\isasymin}\ list{\isacharparenleft}{\kern0pt}A{\isacharparenright}{\kern0pt}{\isacharsemicolon}{\kern0pt}\ P{\isacharparenleft}{\kern0pt}l{\isacharparenright}{\kern0pt}{\isasymrbrakk}\ {\isasymLongrightarrow}\isanewline
\ \ P{\isacharparenleft}{\kern0pt}Cons{\isacharparenleft}{\kern0pt}a{\isacharcomma}{\kern0pt}\ l{\isacharparenright}{\kern0pt}{\isacharparenright}{\kern0pt}{\isasymrbrakk}
{\isasymLongrightarrow}\ P{\isacharparenleft}{\kern0pt}x{\isacharparenright}{\kern0pt}%
\end{isabelle}%
Length, concatenation, and \isa{n}th element of lists are
recursively characterized as follows.%
\begin{isabelle}%
length{\isacharparenleft}{\kern0pt}{\isacharbrackleft}{\kern0pt}{\isacharbrackright}{\kern0pt}{\isacharparenright}{\kern0pt}\ {\isacharequal}{\kern0pt}\ {\isadigit{0}}\isasep\isanewline%
length{\isacharparenleft}{\kern0pt}Cons{\isacharparenleft}{\kern0pt}a{\isacharcomma}{\kern0pt}\ l{\isacharparenright}{\kern0pt}{\isacharparenright}{\kern0pt}\ {\isacharequal}{\kern0pt}\ succ{\isacharparenleft}{\kern0pt}length{\isacharparenleft}{\kern0pt}l{\isacharparenright}{\kern0pt}{\isacharparenright}{\kern0pt}\isasep\isanewline\isanewline%
{\isacharbrackleft}{\kern0pt}{\isacharbrackright}{\kern0pt}\ {\isacharat}{\kern0pt}\ ys\ {\isacharequal}{\kern0pt}\ ys\isasep\isanewline%
Cons{\isacharparenleft}{\kern0pt}a{\isacharcomma}{\kern0pt}\ l{\isacharparenright}{\kern0pt}\ {\isacharat}{\kern0pt}\ ys\ {\isacharequal}{\kern0pt}\ Cons{\isacharparenleft}{\kern0pt}a{\isacharcomma}{\kern0pt}\ l\ {\isacharat}{\kern0pt}\ ys{\isacharparenright}{\kern0pt}\isasep\isanewline\isanewline%
nth{\isacharparenleft}{\kern0pt}{\isadigit{0}}{\isacharcomma}{\kern0pt}\ Cons{\isacharparenleft}{\kern0pt}a{\isacharcomma}{\kern0pt}\ l{\isacharparenright}{\kern0pt}{\isacharparenright}{\kern0pt}\ {\isacharequal}{\kern0pt}\ a\isasep\isanewline%
n\ {\isasymin}\ {\isasymomega}\ {\isasymLongrightarrow}\ nth{\isacharparenleft}{\kern0pt}succ{\isacharparenleft}{\kern0pt}n{\isacharparenright}{\kern0pt}{\isacharcomma}{\kern0pt}\ Cons{\isacharparenleft}{\kern0pt}a{\isacharcomma}{\kern0pt}\ l{\isacharparenright}{\kern0pt}{\isacharparenright}{\kern0pt}\ {\isacharequal}{\kern0pt}\ nth{\isacharparenleft}{\kern0pt}n{\isacharcomma}{\kern0pt}\ l{\isacharparenright}{\kern0pt}%
\end{isabelle}%
We have the usual Haskell-like notation for iterated applications
of \isa{Cons}:%
\begin{isabelle}
\isacommand{lemma}\isamarkupfalse%
\ Cons{\isacharunderscore}{\kern0pt}app{\isacharcolon}{\kern0pt}\ {\isachardoublequoteopen}{\isacharbrackleft}{\kern0pt}a{\isacharcomma}{\kern0pt}b{\isacharcomma}{\kern0pt}c{\isacharbrackright}{\kern0pt}\ {\isacharequal}{\kern0pt}\ Cons{\isacharparenleft}{\kern0pt}a{\isacharcomma}{\kern0pt}Cons{\isacharparenleft}{\kern0pt}b{\isacharcomma}{\kern0pt}Cons{\isacharparenleft}{\kern0pt}c{\isacharcomma}{\kern0pt}{\isacharbrackleft}{\kern0pt}{\isacharbrackright}{\kern0pt}{\isacharparenright}{\kern0pt}{\isacharparenright}{\kern0pt}{\isacharparenright}{\kern0pt}{\isachardoublequoteclose}%
\end{isabelle}

Relative quantifiers restrict the range of the bound variable to a
class \isa{M} of type \isa{i\ {\isasymRightarrow}\ o}; that is, a truth-valued function with
set arguments.%
\begin{isabelle}
\isacommand{lemma}\isamarkupfalse%
\ {\isachardoublequoteopen}{\isasymforall}x{\isacharbrackleft}{\kern0pt}M{\isacharbrackright}{\kern0pt}{\isachardot}{\kern0pt}\ P{\isacharparenleft}{\kern0pt}x{\isacharparenright}{\kern0pt}\ {\isasymequiv}\ {\isasymforall}x{\isachardot}{\kern0pt}\ M{\isacharparenleft}{\kern0pt}x{\isacharparenright}{\kern0pt}\ {\isasymlongrightarrow}\ P{\isacharparenleft}{\kern0pt}x{\isacharparenright}{\kern0pt}{\isachardoublequoteclose}\isanewline
\ \ \ \ \ \ {\isachardoublequoteopen}{\isasymexists}x{\isacharbrackleft}{\kern0pt}M{\isacharbrackright}{\kern0pt}{\isachardot}{\kern0pt}\ P{\isacharparenleft}{\kern0pt}x{\isacharparenright}{\kern0pt}\ {\isasymequiv}\ {\isasymexists}x{\isachardot}{\kern0pt}\ M{\isacharparenleft}{\kern0pt}x{\isacharparenright}{\kern0pt}\ {\isasymand}\ P{\isacharparenleft}{\kern0pt}x{\isacharparenright}{\kern0pt}{\isachardoublequoteclose}
\end{isabelle}
Finally, a set can be viewed (“cast”) as a class using the
following function of type \isa{i\ {\isasymRightarrow}\ i\ {\isasymRightarrow}\ o}.%
\begin{isabelle}%
{\isacharparenleft}{\kern0pt}{\isacharhash}{\kern0pt}{\isacharhash}{\kern0pt}A{\isacharparenright}{\kern0pt}{\isacharparenleft}{\kern0pt}x{\isacharparenright}{\kern0pt}\ {\isasymlongleftrightarrow}\ x\ {\isasymin}\ A%
\end{isabelle}%
\subsection{Relative concepts\label{sec:def-main-relative}%
}
A list of relative concepts (mostly from the \session{ZF-Constructible}
    library) follows next.%
\begin{isabelle}%
big{\isacharunderscore}{\kern0pt}union{\isacharparenleft}{\kern0pt}M{\isacharcomma}{\kern0pt}\ A{\isacharcomma}{\kern0pt}\ z{\isacharparenright}{\kern0pt}\ {\isasymequiv}\ {\isasymforall}x{\isacharbrackleft}{\kern0pt}M{\isacharbrackright}{\kern0pt}{\isachardot}{\kern0pt}\ x\ {\isasymin}\ z\ {\isasymlongleftrightarrow}\ {\isacharparenleft}{\kern0pt}{\isasymexists}y{\isacharbrackleft}{\kern0pt}M{\isacharbrackright}{\kern0pt}{\isachardot}{\kern0pt}\ y\ {\isasymin}\ A\ {\isasymand}\ x\ {\isasymin}\ y{\isacharparenright}{\kern0pt}%
\end{isabelle}%
\begin{isabelle}%
upair{\isacharparenleft}{\kern0pt}M{\isacharcomma}{\kern0pt}\ a{\isacharcomma}{\kern0pt}\ b{\isacharcomma}{\kern0pt}\ z{\isacharparenright}{\kern0pt}\ {\isasymequiv}\ a\ {\isasymin}\ z\ {\isasymand}\ b\ {\isasymin}\ z\ {\isasymand}\ {\isacharparenleft}{\kern0pt}{\isasymforall}x{\isacharbrackleft}{\kern0pt}M{\isacharbrackright}{\kern0pt}{\isachardot}{\kern0pt}\ x\ {\isasymin}\ z\ {\isasymlongrightarrow}\ x\ {\isacharequal}{\kern0pt}\ a\ {\isasymor}\ x\ {\isacharequal}{\kern0pt}\ b{\isacharparenright}{\kern0pt}%
\end{isabelle}%
\begin{isabelle}%
pair{\isacharparenleft}{\kern0pt}M{\isacharcomma}{\kern0pt}\ a{\isacharcomma}{\kern0pt}\ b{\isacharcomma}{\kern0pt}\ z{\isacharparenright}{\kern0pt}\ {\isasymequiv}\ 
{\isasymexists}x{\isacharbrackleft}{\kern0pt}M{\isacharbrackright}{\kern0pt}{\isachardot}{\kern0pt}\ upair{\isacharparenleft}{\kern0pt}M{\isacharcomma}{\kern0pt}\ a{\isacharcomma}{\kern0pt}\ a{\isacharcomma}{\kern0pt}\ x{\isacharparenright}{\kern0pt}\ {\isasymand}\isanewline
\ \ \ \ \ \ \ \ \ \ \ \ \ \ \ \ \ \ \ \ \ \ {\isacharparenleft}{\kern0pt}{\isasymexists}y{\isacharbrackleft}{\kern0pt}M{\isacharbrackright}{\kern0pt}{\isachardot}{\kern0pt}\ upair{\isacharparenleft}{\kern0pt}M{\isacharcomma}{\kern0pt}\ a{\isacharcomma}{\kern0pt}\ b{\isacharcomma}{\kern0pt}\ y{\isacharparenright}{\kern0pt}\ {\isasymand}\ upair{\isacharparenleft}{\kern0pt}M{\isacharcomma}{\kern0pt}\ x{\isacharcomma}{\kern0pt}\ y{\isacharcomma}{\kern0pt}\ z{\isacharparenright}{\kern0pt}{\isacharparenright}{\kern0pt}%
\end{isabelle}%
\begin{isabelle}%
successor{\isacharparenleft}{\kern0pt}M{\isacharcomma}{\kern0pt}\ a{\isacharcomma}{\kern0pt}\ z{\isacharparenright}{\kern0pt}\ {\isasymequiv}\isanewline
{\isasymexists}x{\isacharbrackleft}{\kern0pt}M{\isacharbrackright}{\kern0pt}{\isachardot}{\kern0pt}\ upair{\isacharparenleft}{\kern0pt}M{\isacharcomma}{\kern0pt}\ a{\isacharcomma}{\kern0pt}\ a{\isacharcomma}{\kern0pt}\ x{\isacharparenright}{\kern0pt}\ {\isasymand}\ {\isacharparenleft}{\kern0pt}{\isasymforall}xa{\isacharbrackleft}{\kern0pt}M{\isacharbrackright}{\kern0pt}{\isachardot}{\kern0pt}\ xa\ {\isasymin}\ z\ {\isasymlongleftrightarrow}\ xa\ {\isasymin}\ x\ {\isasymor}\ xa\ {\isasymin}\ a{\isacharparenright}{\kern0pt}%
\end{isabelle}%
\begin{isabelle}%
empty{\isacharparenleft}{\kern0pt}M{\isacharcomma}{\kern0pt}\ z{\isacharparenright}{\kern0pt}\ {\isasymequiv}\ {\isasymforall}x{\isacharbrackleft}{\kern0pt}M{\isacharbrackright}{\kern0pt}{\isachardot}{\kern0pt}\ x\ {\isasymnotin}\ z%
\end{isabelle}%
\begin{isabelle}%
transitive{\isacharunderscore}{\kern0pt}set{\isacharparenleft}{\kern0pt}M{\isacharcomma}{\kern0pt}\ a{\isacharparenright}{\kern0pt}\ {\isasymequiv}\ {\isasymforall}x{\isacharbrackleft}{\kern0pt}M{\isacharbrackright}{\kern0pt}{\isachardot}{\kern0pt}\ x\ {\isasymin}\ a\ {\isasymlongrightarrow}\ {\isacharparenleft}{\kern0pt}{\isasymforall}xa{\isacharbrackleft}{\kern0pt}M{\isacharbrackright}{\kern0pt}{\isachardot}{\kern0pt}\ xa\ {\isasymin}\ x\ {\isasymlongrightarrow}\ xa\ {\isasymin}\ a{\isacharparenright}{\kern0pt}%
\end{isabelle}%
\begin{isabelle}%
ordinal{\isacharparenleft}{\kern0pt}M{\isacharcomma}{\kern0pt}\ a{\isacharparenright}{\kern0pt}\ {\isasymequiv}\isanewline
transitive{\isacharunderscore}{\kern0pt}set{\isacharparenleft}{\kern0pt}M{\isacharcomma}{\kern0pt}\ a{\isacharparenright}{\kern0pt}\ {\isasymand}\ {\isacharparenleft}{\kern0pt}{\isasymforall}x{\isacharbrackleft}{\kern0pt}M{\isacharbrackright}{\kern0pt}{\isachardot}{\kern0pt}\ x\ {\isasymin}\ a\ {\isasymlongrightarrow}\ transitive{\isacharunderscore}{\kern0pt}set{\isacharparenleft}{\kern0pt}M{\isacharcomma}{\kern0pt}\ x{\isacharparenright}{\kern0pt}{\isacharparenright}{\kern0pt}%
\end{isabelle}%
\begin{isabelle}%
image{\isacharparenleft}{\kern0pt}M{\isacharcomma}{\kern0pt}\ r{\isacharcomma}{\kern0pt}\ A{\isacharcomma}{\kern0pt}\ z{\isacharparenright}{\kern0pt}\ {\isasymequiv}\isanewline
{\isasymforall}y{\isacharbrackleft}{\kern0pt}M{\isacharbrackright}{\kern0pt}{\isachardot}{\kern0pt}\ y\ {\isasymin}\ z\ {\isasymlongleftrightarrow}\ {\isacharparenleft}{\kern0pt}{\isasymexists}w{\isacharbrackleft}{\kern0pt}M{\isacharbrackright}{\kern0pt}{\isachardot}{\kern0pt}\ w\ {\isasymin}\ r\ {\isasymand}\ {\isacharparenleft}{\kern0pt}{\isasymexists}x{\isacharbrackleft}{\kern0pt}M{\isacharbrackright}{\kern0pt}{\isachardot}{\kern0pt}\ x\ {\isasymin}\ A\ {\isasymand}\ pair{\isacharparenleft}{\kern0pt}M{\isacharcomma}{\kern0pt}\ x{\isacharcomma}{\kern0pt}\ y{\isacharcomma}{\kern0pt}\ w{\isacharparenright}{\kern0pt}{\isacharparenright}{\kern0pt}{\isacharparenright}{\kern0pt}%
\end{isabelle}%
\begin{isabelle}%
is{\isacharunderscore}{\kern0pt}apply{\isacharparenleft}{\kern0pt}M{\isacharcomma}{\kern0pt}\ f{\isacharcomma}{\kern0pt}\ x{\isacharcomma}{\kern0pt}\ y{\isacharparenright}{\kern0pt}\ {\isasymequiv}\isanewline
{\isasymexists}xs{\isacharbrackleft}{\kern0pt}M{\isacharbrackright}{\kern0pt}{\isachardot}{\kern0pt}\isanewline
\isaindent{\ \ \ }{\isasymexists}fxs{\isacharbrackleft}{\kern0pt}M{\isacharbrackright}{\kern0pt}{\isachardot}{\kern0pt}\ upair{\isacharparenleft}{\kern0pt}M{\isacharcomma}{\kern0pt}\ x{\isacharcomma}{\kern0pt}\ x{\isacharcomma}{\kern0pt}\ xs{\isacharparenright}{\kern0pt}\ {\isasymand}\ image{\isacharparenleft}{\kern0pt}M{\isacharcomma}{\kern0pt}\ f{\isacharcomma}{\kern0pt}\ xs{\isacharcomma}{\kern0pt}\ fxs{\isacharparenright}{\kern0pt}\ {\isasymand}\isanewline
\isaindent{\ \ \ \ \ }big{\isacharunderscore}{\kern0pt}union{\isacharparenleft}{\kern0pt}M{\isacharcomma}{\kern0pt}\ fxs{\isacharcomma}{\kern0pt}\ y{\isacharparenright}{\kern0pt}%
\end{isabelle}%
\begin{isabelle}%
is{\isacharunderscore}{\kern0pt}function{\isacharparenleft}{\kern0pt}M{\isacharcomma}{\kern0pt}\ r{\isacharparenright}{\kern0pt}\ {\isasymequiv}\isanewline
{\isasymforall}x{\isacharbrackleft}{\kern0pt}M{\isacharbrackright}{\kern0pt}{\isachardot}{\kern0pt}\isanewline
\isaindent{\ \ \ }{\isasymforall}y{\isacharbrackleft}{\kern0pt}M{\isacharbrackright}{\kern0pt}{\isachardot}{\kern0pt}\isanewline
\isaindent{\ \ \ \ \ \ }{\isasymforall}y{\isacharprime}{\kern0pt}{\isacharbrackleft}{\kern0pt}M{\isacharbrackright}{\kern0pt}{\isachardot}{\kern0pt}\isanewline
\isaindent{\ \ \ \ \ \ \ \ \ }{\isasymforall}p{\isacharbrackleft}{\kern0pt}M{\isacharbrackright}{\kern0pt}{\isachardot}{\kern0pt}\isanewline
\isaindent{\ \ \ \ \ \ \ \ \ \ \ \ }{\isasymforall}p{\isacharprime}{\kern0pt}{\isacharbrackleft}{\kern0pt}M{\isacharbrackright}{\kern0pt}{\isachardot}{\kern0pt}\isanewline
\isaindent{\ \ \ \ \ \ \ \ \ \ \ \ \ \ \ }pair{\isacharparenleft}{\kern0pt}M{\isacharcomma}{\kern0pt}\ x{\isacharcomma}{\kern0pt}\ y{\isacharcomma}{\kern0pt}\ p{\isacharparenright}{\kern0pt}\ {\isasymlongrightarrow}\isanewline
\isaindent{\ \ \ \ \ \ \ \ \ \ \ \ \ \ \ }pair{\isacharparenleft}{\kern0pt}M{\isacharcomma}{\kern0pt}\ x{\isacharcomma}{\kern0pt}\ y{\isacharprime}{\kern0pt}{\isacharcomma}{\kern0pt}\ p{\isacharprime}{\kern0pt}{\isacharparenright}{\kern0pt}\ {\isasymlongrightarrow}\ p\ {\isasymin}\ r\ {\isasymlongrightarrow}\ p{\isacharprime}{\kern0pt}\ {\isasymin}\ r\ {\isasymlongrightarrow}\ y\ {\isacharequal}{\kern0pt}\ y{\isacharprime}{\kern0pt}%
\end{isabelle}%
\begin{isabelle}%
is{\isacharunderscore}{\kern0pt}relation{\isacharparenleft}{\kern0pt}M{\isacharcomma}{\kern0pt}\ r{\isacharparenright}{\kern0pt}\ {\isasymequiv}\ {\isasymforall}z{\isacharbrackleft}{\kern0pt}M{\isacharbrackright}{\kern0pt}{\isachardot}{\kern0pt}\ z\ {\isasymin}\ r\ {\isasymlongrightarrow}\ {\isacharparenleft}{\kern0pt}{\isasymexists}x{\isacharbrackleft}{\kern0pt}M{\isacharbrackright}{\kern0pt}{\isachardot}{\kern0pt}\ {\isasymexists}y{\isacharbrackleft}{\kern0pt}M{\isacharbrackright}{\kern0pt}{\isachardot}{\kern0pt}\ pair{\isacharparenleft}{\kern0pt}M{\isacharcomma}{\kern0pt}\ x{\isacharcomma}{\kern0pt}\ y{\isacharcomma}{\kern0pt}\ z{\isacharparenright}{\kern0pt}{\isacharparenright}{\kern0pt}%
\end{isabelle}%
\begin{isabelle}%
is{\isacharunderscore}{\kern0pt}domain{\isacharparenleft}{\kern0pt}M{\isacharcomma}{\kern0pt}\ r{\isacharcomma}{\kern0pt}\ z{\isacharparenright}{\kern0pt}\ {\isasymequiv}\isanewline
{\isasymforall}x{\isacharbrackleft}{\kern0pt}M{\isacharbrackright}{\kern0pt}{\isachardot}{\kern0pt}\ x\ {\isasymin}\ z\ {\isasymlongleftrightarrow}\ {\isacharparenleft}{\kern0pt}{\isasymexists}w{\isacharbrackleft}{\kern0pt}M{\isacharbrackright}{\kern0pt}{\isachardot}{\kern0pt}\ w\ {\isasymin}\ r\ {\isasymand}\ {\isacharparenleft}{\kern0pt}{\isasymexists}y{\isacharbrackleft}{\kern0pt}M{\isacharbrackright}{\kern0pt}{\isachardot}{\kern0pt}\ pair{\isacharparenleft}{\kern0pt}M{\isacharcomma}{\kern0pt}\ x{\isacharcomma}{\kern0pt}\ y{\isacharcomma}{\kern0pt}\ w{\isacharparenright}{\kern0pt}{\isacharparenright}{\kern0pt}{\isacharparenright}{\kern0pt}%
\end{isabelle}%
\begin{isabelle}%
typed{\isacharunderscore}{\kern0pt}function{\isacharparenleft}{\kern0pt}M{\isacharcomma}{\kern0pt}\ A{\isacharcomma}{\kern0pt}\ B{\isacharcomma}{\kern0pt}\ r{\isacharparenright}{\kern0pt}\ {\isasymequiv}\isanewline
is{\isacharunderscore}{\kern0pt}function{\isacharparenleft}{\kern0pt}M{\isacharcomma}{\kern0pt}\ r{\isacharparenright}{\kern0pt}\ {\isasymand}\isanewline
is{\isacharunderscore}{\kern0pt}relation{\isacharparenleft}{\kern0pt}M{\isacharcomma}{\kern0pt}\ r{\isacharparenright}{\kern0pt}\ {\isasymand}\isanewline
is{\isacharunderscore}{\kern0pt}domain{\isacharparenleft}{\kern0pt}M{\isacharcomma}{\kern0pt}\ r{\isacharcomma}{\kern0pt}\ A{\isacharparenright}{\kern0pt}\ {\isasymand}\isanewline
{\isacharparenleft}{\kern0pt}{\isasymforall}u{\isacharbrackleft}{\kern0pt}M{\isacharbrackright}{\kern0pt}{\isachardot}{\kern0pt}\ u\ {\isasymin}\ r\ {\isasymlongrightarrow}\ {\isacharparenleft}{\kern0pt}{\isasymforall}x{\isacharbrackleft}{\kern0pt}M{\isacharbrackright}{\kern0pt}{\isachardot}{\kern0pt}\ {\isasymforall}y{\isacharbrackleft}{\kern0pt}M{\isacharbrackright}{\kern0pt}{\isachardot}{\kern0pt}\ pair{\isacharparenleft}{\kern0pt}M{\isacharcomma}{\kern0pt}\ x{\isacharcomma}{\kern0pt}\ y{\isacharcomma}{\kern0pt}\ u{\isacharparenright}{\kern0pt}\ {\isasymlongrightarrow}\ y\ {\isasymin}\ B{\isacharparenright}{\kern0pt}{\isacharparenright}{\kern0pt}%
\end{isabelle}%
\begin{isabelle}%
is{\isacharunderscore}{\kern0pt}function{\isacharunderscore}{\kern0pt}space{\isacharparenleft}{\kern0pt}M{\isacharcomma}{\kern0pt}\ A{\isacharcomma}{\kern0pt}\ B{\isacharcomma}{\kern0pt}\ fs{\isacharparenright}{\kern0pt}\ {\isasymequiv}\isanewline
M{\isacharparenleft}{\kern0pt}fs{\isacharparenright}{\kern0pt}\ {\isasymand}\ {\isacharparenleft}{\kern0pt}{\isasymforall}f{\isacharbrackleft}{\kern0pt}M{\isacharbrackright}{\kern0pt}{\isachardot}{\kern0pt}\ f\ {\isasymin}\ fs\ {\isasymlongleftrightarrow}\ typed{\isacharunderscore}{\kern0pt}function{\isacharparenleft}{\kern0pt}M{\isacharcomma}{\kern0pt}\ A{\isacharcomma}{\kern0pt}\ B{\isacharcomma}{\kern0pt}\ f{\isacharparenright}{\kern0pt}{\isacharparenright}{\kern0pt}\isasep\isanewline\isanewline%
A\ {\isasymrightarrow}\isactrlbsup M\isactrlesup \ B\ {\isasymequiv}\ THE\ d{\isachardot}{\kern0pt}\ is{\isacharunderscore}{\kern0pt}function{\isacharunderscore}{\kern0pt}space{\isacharparenleft}{\kern0pt}M{\isacharcomma}{\kern0pt}\ A{\isacharcomma}{\kern0pt}\ B{\isacharcomma}{\kern0pt}\ d{\isacharparenright}{\kern0pt}\isasep\isanewline\isanewline%
surjection{\isacharparenleft}{\kern0pt}M{\isacharcomma}{\kern0pt}\ A{\isacharcomma}{\kern0pt}\ B{\isacharcomma}{\kern0pt}\ f{\isacharparenright}{\kern0pt}\ {\isasymequiv}\isanewline
typed{\isacharunderscore}{\kern0pt}function{\isacharparenleft}{\kern0pt}M{\isacharcomma}{\kern0pt}\ A{\isacharcomma}{\kern0pt}\ B{\isacharcomma}{\kern0pt}\ f{\isacharparenright}{\kern0pt}\ {\isasymand}\isanewline
{\isacharparenleft}{\kern0pt}{\isasymforall}y{\isacharbrackleft}{\kern0pt}M{\isacharbrackright}{\kern0pt}{\isachardot}{\kern0pt}\ y\ {\isasymin}\ B\ {\isasymlongrightarrow}\ {\isacharparenleft}{\kern0pt}{\isasymexists}x{\isacharbrackleft}{\kern0pt}M{\isacharbrackright}{\kern0pt}{\isachardot}{\kern0pt}\ x\ {\isasymin}\ A\ {\isasymand}\ is{\isacharunderscore}{\kern0pt}apply{\isacharparenleft}{\kern0pt}M{\isacharcomma}{\kern0pt}\ f{\isacharcomma}{\kern0pt}\ x{\isacharcomma}{\kern0pt}\ y{\isacharparenright}{\kern0pt}{\isacharparenright}{\kern0pt}{\isacharparenright}{\kern0pt}%
\end{isabelle}%

\subsubsection*{Relative version of the $\ZFC$ axioms}
\begin{isabelle}%
extensionality{\isacharparenleft}{\kern0pt}M{\isacharparenright}{\kern0pt}\ {\isasymequiv}\ {\isasymforall}x{\isacharbrackleft}{\kern0pt}M{\isacharbrackright}{\kern0pt}{\isachardot}{\kern0pt}\ {\isasymforall}y{\isacharbrackleft}{\kern0pt}M{\isacharbrackright}{\kern0pt}{\isachardot}{\kern0pt}\ {\isacharparenleft}{\kern0pt}{\isasymforall}z{\isacharbrackleft}{\kern0pt}M{\isacharbrackright}{\kern0pt}{\isachardot}{\kern0pt}\ z\ {\isasymin}\ x\ {\isasymlongleftrightarrow}\ z\ {\isasymin}\ y{\isacharparenright}{\kern0pt}\ {\isasymlongrightarrow}\ x\ {\isacharequal}{\kern0pt}\ y%
\end{isabelle}%
\begin{isabelle}%
foundation{\isacharunderscore}{\kern0pt}ax{\isacharparenleft}{\kern0pt}M{\isacharparenright}{\kern0pt}\ {\isasymequiv}\isanewline
{\isasymforall}x{\isacharbrackleft}{\kern0pt}M{\isacharbrackright}{\kern0pt}{\isachardot}{\kern0pt}\ {\isacharparenleft}{\kern0pt}{\isasymexists}y{\isacharbrackleft}{\kern0pt}M{\isacharbrackright}{\kern0pt}{\isachardot}{\kern0pt}\ y\ {\isasymin}\ x{\isacharparenright}{\kern0pt}\ {\isasymlongrightarrow}\ {\isacharparenleft}{\kern0pt}{\isasymexists}y{\isacharbrackleft}{\kern0pt}M{\isacharbrackright}{\kern0pt}{\isachardot}{\kern0pt}\ y\ {\isasymin}\ x\ {\isasymand}\ {\isasymnot}\ {\isacharparenleft}{\kern0pt}{\isasymexists}z{\isacharbrackleft}{\kern0pt}M{\isacharbrackright}{\kern0pt}{\isachardot}{\kern0pt}\ z\ {\isasymin}\ x\ {\isasymand}\ z\ {\isasymin}\ y{\isacharparenright}{\kern0pt}{\isacharparenright}{\kern0pt}%
\end{isabelle}%
\begin{isabelle}%
upair{\isacharunderscore}{\kern0pt}ax{\isacharparenleft}{\kern0pt}M{\isacharparenright}{\kern0pt}\ {\isasymequiv}\ {\isasymforall}x{\isacharbrackleft}{\kern0pt}M{\isacharbrackright}{\kern0pt}{\isachardot}{\kern0pt}\ {\isasymforall}y{\isacharbrackleft}{\kern0pt}M{\isacharbrackright}{\kern0pt}{\isachardot}{\kern0pt}\ {\isasymexists}z{\isacharbrackleft}{\kern0pt}M{\isacharbrackright}{\kern0pt}{\isachardot}{\kern0pt}\ upair{\isacharparenleft}{\kern0pt}M{\isacharcomma}{\kern0pt}\ x{\isacharcomma}{\kern0pt}\ y{\isacharcomma}{\kern0pt}\ z{\isacharparenright}{\kern0pt}%
\end{isabelle}%
\begin{isabelle}%
Union{\isacharunderscore}{\kern0pt}ax{\isacharparenleft}{\kern0pt}M{\isacharparenright}{\kern0pt}\ {\isasymequiv}\ {\isasymforall}x{\isacharbrackleft}{\kern0pt}M{\isacharbrackright}{\kern0pt}{\isachardot}{\kern0pt}\ {\isasymexists}z{\isacharbrackleft}{\kern0pt}M{\isacharbrackright}{\kern0pt}{\isachardot}{\kern0pt}\ big{\isacharunderscore}{\kern0pt}union{\isacharparenleft}{\kern0pt}M{\isacharcomma}{\kern0pt}\ x{\isacharcomma}{\kern0pt}\ z{\isacharparenright}{\kern0pt}%
\end{isabelle}%
\begin{isabelle}%
power{\isacharunderscore}{\kern0pt}ax{\isacharparenleft}{\kern0pt}M{\isacharparenright}{\kern0pt}\ {\isasymequiv}\ {\isasymforall}x{\isacharbrackleft}{\kern0pt}M{\isacharbrackright}{\kern0pt}{\isachardot}{\kern0pt}\ {\isasymexists}z{\isacharbrackleft}{\kern0pt}M{\isacharbrackright}{\kern0pt}{\isachardot}{\kern0pt}\ {\isasymforall}xa{\isacharbrackleft}{\kern0pt}M{\isacharbrackright}{\kern0pt}{\isachardot}{\kern0pt}\ xa\ {\isasymin}\ z\ {\isasymlongleftrightarrow}\isanewline
\ \ \ \ \ \ \ \ \ \ \ \ \ \ \ \ \ {\isacharparenleft}{\kern0pt}{\isasymforall}xb{\isacharbrackleft}{\kern0pt}M{\isacharbrackright}{\kern0pt}{\isachardot}{\kern0pt}\ xb\ {\isasymin}\ xa\ {\isasymlongrightarrow}\ xb\ {\isasymin}\ x{\isacharparenright}{\kern0pt}%
\end{isabelle}%
\begin{isabelle}%
infinity{\isacharunderscore}{\kern0pt}ax{\isacharparenleft}{\kern0pt}M{\isacharparenright}{\kern0pt}\ {\isasymequiv}\isanewline
{\isasymexists}I{\isacharbrackleft}{\kern0pt}M{\isacharbrackright}{\kern0pt}{\isachardot}{\kern0pt}\isanewline
\isaindent{\ \ \ }{\isacharparenleft}{\kern0pt}{\isasymexists}z{\isacharbrackleft}{\kern0pt}M{\isacharbrackright}{\kern0pt}{\isachardot}{\kern0pt}\ empty{\isacharparenleft}{\kern0pt}M{\isacharcomma}{\kern0pt}\ z{\isacharparenright}{\kern0pt}\ {\isasymand}\ z\ {\isasymin}\ I{\isacharparenright}{\kern0pt}\ {\isasymand}\isanewline
\isaindent{\ \ \ }{\isacharparenleft}{\kern0pt}{\isasymforall}y{\isacharbrackleft}{\kern0pt}M{\isacharbrackright}{\kern0pt}{\isachardot}{\kern0pt}\ y\ {\isasymin}\ I\ {\isasymlongrightarrow}\ {\isacharparenleft}{\kern0pt}{\isasymexists}sy{\isacharbrackleft}{\kern0pt}M{\isacharbrackright}{\kern0pt}{\isachardot}{\kern0pt}\ successor{\isacharparenleft}{\kern0pt}M{\isacharcomma}{\kern0pt}\ y{\isacharcomma}{\kern0pt}\ sy{\isacharparenright}{\kern0pt}\ {\isasymand}\ sy\ {\isasymin}\ I{\isacharparenright}{\kern0pt}{\isacharparenright}{\kern0pt}%
\end{isabelle}%
\begin{isabelle}%
choice{\isacharunderscore}{\kern0pt}ax{\isacharparenleft}{\kern0pt}M{\isacharparenright}{\kern0pt}\ {\isasymequiv}\ {\isasymforall}x{\isacharbrackleft}{\kern0pt}M{\isacharbrackright}{\kern0pt}{\isachardot}{\kern0pt}\ {\isasymexists}a{\isacharbrackleft}{\kern0pt}M{\isacharbrackright}{\kern0pt}{\isachardot}{\kern0pt}\ {\isasymexists}f{\isacharbrackleft}{\kern0pt}M{\isacharbrackright}{\kern0pt}{\isachardot}{\kern0pt}\ ordinal{\isacharparenleft}{\kern0pt}M{\isacharcomma}{\kern0pt}\ a{\isacharparenright}{\kern0pt}\ {\isasymand}\isanewline
\ \ \ \ \ \ \ \ \ \ \ \ \ \ \ \ \ \ surjection{\isacharparenleft}{\kern0pt}M{\isacharcomma}{\kern0pt}\ a{\isacharcomma}{\kern0pt}\ x{\isacharcomma}{\kern0pt}\ f{\isacharparenright}{\kern0pt}%
\end{isabelle}%
\begin{isabelle}%
separation{\isacharparenleft}{\kern0pt}M{\isacharcomma}{\kern0pt}\ P{\isacharparenright}{\kern0pt}\ {\isasymequiv}\ {\isasymforall}z{\isacharbrackleft}{\kern0pt}M{\isacharbrackright}{\kern0pt}{\isachardot}{\kern0pt}\ {\isasymexists}y{\isacharbrackleft}{\kern0pt}M{\isacharbrackright}{\kern0pt}{\isachardot}{\kern0pt}\ {\isasymforall}x{\isacharbrackleft}{\kern0pt}M{\isacharbrackright}{\kern0pt}{\isachardot}{\kern0pt}\ x\ {\isasymin}\ y\ {\isasymlongleftrightarrow}\ x\ {\isasymin}\ z\ {\isasymand}\ P{\isacharparenleft}{\kern0pt}x{\isacharparenright}{\kern0pt}%
\end{isabelle}%
\begin{isabelle}%
univalent{\isacharparenleft}{\kern0pt}M{\isacharcomma}{\kern0pt}\ A{\isacharcomma}{\kern0pt}\ P{\isacharparenright}{\kern0pt}\ {\isasymequiv}\isanewline
{\isasymforall}x{\isacharbrackleft}{\kern0pt}M{\isacharbrackright}{\kern0pt}{\isachardot}{\kern0pt}\ x\ {\isasymin}\ A\ {\isasymlongrightarrow}\ {\isacharparenleft}{\kern0pt}{\isasymforall}y{\isacharbrackleft}{\kern0pt}M{\isacharbrackright}{\kern0pt}{\isachardot}{\kern0pt}\ {\isasymforall}z{\isacharbrackleft}{\kern0pt}M{\isacharbrackright}{\kern0pt}{\isachardot}{\kern0pt}\ P{\isacharparenleft}{\kern0pt}x{\isacharcomma}{\kern0pt}\ y{\isacharparenright}{\kern0pt}\ {\isasymand}\ P{\isacharparenleft}{\kern0pt}x{\isacharcomma}{\kern0pt}\ z{\isacharparenright}{\kern0pt}\ {\isasymlongrightarrow}\ y\ {\isacharequal}{\kern0pt}\ z{\isacharparenright}{\kern0pt}%
\end{isabelle}%
\begin{isabelle}%
strong{\isacharunderscore}{\kern0pt}replacement{\isacharparenleft}{\kern0pt}M{\isacharcomma}{\kern0pt}\ P{\isacharparenright}{\kern0pt}\ {\isasymequiv}\isanewline
{\isasymforall}A{\isacharbrackleft}{\kern0pt}M{\isacharbrackright}{\kern0pt}{\isachardot}{\kern0pt}\isanewline
\isaindent{\ \ \ }univalent{\isacharparenleft}{\kern0pt}M{\isacharcomma}{\kern0pt}\ A{\isacharcomma}{\kern0pt}\ P{\isacharparenright}{\kern0pt}\ {\isasymlongrightarrow}\ {\isacharparenleft}{\kern0pt}{\isasymexists}Y{\isacharbrackleft}{\kern0pt}M{\isacharbrackright}{\kern0pt}{\isachardot}{\kern0pt}\ {\isasymforall}b{\isacharbrackleft}{\kern0pt}M{\isacharbrackright}{\kern0pt}{\isachardot}{\kern0pt}\ b\ {\isasymin}\ Y\ {\isasymlongleftrightarrow}\isanewline
\ \ \ \ \ \ {\isacharparenleft}{\kern0pt}{\isasymexists}x{\isacharbrackleft}{\kern0pt}M{\isacharbrackright}{\kern0pt}{\isachardot}{\kern0pt}\ x\ {\isasymin}\ A\ {\isasymand}
P{\isacharparenleft}{\kern0pt}x{\isacharcomma}{\kern0pt}\ b{\isacharparenright}{\kern0pt}{\isacharparenright}{\kern0pt}{\isacharparenright}{\kern0pt}%
\end{isabelle}%
\subsubsection*{Internalized formulas}
“Codes” for formulas (as sets) are constructed from natural
numbers using \isa{Member}, \isa{Equal}, \isa{Nand},
and \isa{Forall}.%
\begin{isabelle}%
{\isasymlbrakk}x\ {\isasymin}\ {\isasymomega}{\isacharsemicolon}{\kern0pt}\ y\ {\isasymin}\ {\isasymomega}{\isasymrbrakk}\ {\isasymLongrightarrow}\ {\isasymcdot}x\ {\isasymin}\ y{\isasymcdot}\ {\isasymin}\ formula\isasep\isanewline%
{\isasymlbrakk}x\ {\isasymin}\ {\isasymomega}{\isacharsemicolon}{\kern0pt}\ y\ {\isasymin}\ {\isasymomega}{\isasymrbrakk}\ {\isasymLongrightarrow}\ {\isasymcdot}x\ {\isacharequal}{\kern0pt}\ y{\isasymcdot}\ {\isasymin}\ formula\isasep\isanewline%
{\isasymlbrakk}p\ {\isasymin}\ formula{\isacharsemicolon}{\kern0pt}\ q\ {\isasymin}\ formula{\isasymrbrakk}\ {\isasymLongrightarrow}\ {\isasymcdot}{\isasymnot}{\isacharparenleft}{\kern0pt}p\ {\isasymand}\ q{\isacharparenright}{\kern0pt}{\isasymcdot}\ {\isasymin}\ formula\isasep\isanewline%
p\ {\isasymin}\ formula\ {\isasymLongrightarrow}\ {\isacharparenleft}{\kern0pt}{\isasymcdot}{\isasymforall}p{\isasymcdot}{\isacharparenright}{\kern0pt}\ {\isasymin}\ formula\isasep\isanewline\isanewline%
{\isasymlbrakk}x\ {\isasymin}\ formula{\isacharsemicolon}{\kern0pt}\ {\isasymAnd}x\ y{\isachardot}{\kern0pt}\ {\isasymlbrakk}x\ {\isasymin}\ {\isasymomega}{\isacharsemicolon}{\kern0pt}\ y\ {\isasymin}\ {\isasymomega}{\isasymrbrakk}\ {\isasymLongrightarrow}\ P{\isacharparenleft}{\kern0pt}{\isasymcdot}x\ {\isasymin}\ y{\isasymcdot}{\isacharparenright}{\kern0pt}{\isacharsemicolon}{\kern0pt}\isanewline
\isaindent{\ }{\isasymAnd}x\ y{\isachardot}{\kern0pt}\ {\isasymlbrakk}x\ {\isasymin}\ {\isasymomega}{\isacharsemicolon}{\kern0pt}\ y\ {\isasymin}\ {\isasymomega}{\isasymrbrakk}\ {\isasymLongrightarrow}\ P{\isacharparenleft}{\kern0pt}{\isasymcdot}x\ {\isacharequal}{\kern0pt}\ y{\isasymcdot}{\isacharparenright}{\kern0pt}{\isacharsemicolon}{\kern0pt}\isanewline
\isaindent{\ }{\isasymAnd}p\ q{\isachardot}{\kern0pt}\ {\isasymlbrakk}p\ {\isasymin}\ formula{\isacharsemicolon}{\kern0pt}\ P{\isacharparenleft}{\kern0pt}p{\isacharparenright}{\kern0pt}{\isacharsemicolon}{\kern0pt}\ q\ {\isasymin}\ formula{\isacharsemicolon}{\kern0pt}\ P{\isacharparenleft}{\kern0pt}q{\isacharparenright}{\kern0pt}{\isasymrbrakk}\ {\isasymLongrightarrow}\ P{\isacharparenleft}{\kern0pt}{\isasymcdot}{\isasymnot}{\isacharparenleft}{\kern0pt}p\ {\isasymand}\ q{\isacharparenright}{\kern0pt}{\isasymcdot}{\isacharparenright}{\kern0pt}{\isacharsemicolon}{\kern0pt}\isanewline
\isaindent{\ }{\isasymAnd}p{\isachardot}{\kern0pt}\ {\isasymlbrakk}p\ {\isasymin}\ formula{\isacharsemicolon}{\kern0pt}\ P{\isacharparenleft}{\kern0pt}p{\isacharparenright}{\kern0pt}{\isasymrbrakk}\ {\isasymLongrightarrow}\ P{\isacharparenleft}{\kern0pt}{\isacharparenleft}{\kern0pt}{\isasymcdot}{\isasymforall}p{\isasymcdot}{\isacharparenright}{\kern0pt}{\isacharparenright}{\kern0pt}{\isasymrbrakk}\isanewline
{\isasymLongrightarrow}\ P{\isacharparenleft}{\kern0pt}x{\isacharparenright}{\kern0pt}%
\end{isabelle}%
Definitions for the other connectives and the internal existential
quantifier are also provided. For instance, negation:%
\begin{isabelle}%
{\isasymcdot}{\isasymnot}p{\isasymcdot}\ {\isasymequiv}\ {\isasymcdot}{\isasymnot}{\isacharparenleft}{\kern0pt}p\ {\isasymand}\ p{\isacharparenright}{\kern0pt}{\isasymcdot}%
\end{isabelle}%
The \isa{arity} function strictly bounding the free de Bruijn
indices of a formula is defined below:
\begin{isabelle}%
arity{\isacharparenleft}{\kern0pt}{\isasymcdot}x\ {\isasymin}\ y{\isasymcdot}{\isacharparenright}{\kern0pt}\ {\isacharequal}{\kern0pt}\ succ{\isacharparenleft}{\kern0pt}x{\isacharparenright}{\kern0pt}\ {\isasymunion}\ succ{\isacharparenleft}{\kern0pt}y{\isacharparenright}{\kern0pt}\isasep\isanewline%
arity{\isacharparenleft}{\kern0pt}{\isasymcdot}x\ {\isacharequal}{\kern0pt}\ y{\isasymcdot}{\isacharparenright}{\kern0pt}\ {\isacharequal}{\kern0pt}\ succ{\isacharparenleft}{\kern0pt}x{\isacharparenright}{\kern0pt}\ {\isasymunion}\ succ{\isacharparenleft}{\kern0pt}y{\isacharparenright}{\kern0pt}\isasep\isanewline%
arity{\isacharparenleft}{\kern0pt}{\isasymcdot}{\isasymnot}{\isacharparenleft}{\kern0pt}p\ {\isasymand}\ q{\isacharparenright}{\kern0pt}{\isasymcdot}{\isacharparenright}{\kern0pt}\ {\isacharequal}{\kern0pt}\ arity{\isacharparenleft}{\kern0pt}p{\isacharparenright}{\kern0pt}\ {\isasymunion}\ arity{\isacharparenleft}{\kern0pt}q{\isacharparenright}{\kern0pt}\isasep\isanewline%
arity{\isacharparenleft}{\kern0pt}{\isacharparenleft}{\kern0pt}{\isasymcdot}{\isasymforall}p{\isasymcdot}{\isacharparenright}{\kern0pt}{\isacharparenright}{\kern0pt}\ {\isacharequal}{\kern0pt}\ pred{\isacharparenleft}{\kern0pt}arity{\isacharparenleft}{\kern0pt}p{\isacharparenright}{\kern0pt}{\isacharparenright}{\kern0pt}%
\end{isabelle}%
We have the satisfaction relation between $\in$-models and
    first order formulas (given a “environment” list representing
    the assignment of free variables),%
\begin{isabelle}%
{\isasymlbrakk}nth{\isacharparenleft}{\kern0pt}i{\isacharcomma}{\kern0pt}\ env{\isacharparenright}{\kern0pt}\ {\isacharequal}{\kern0pt}\ x{\isacharsemicolon}{\kern0pt}\ nth{\isacharparenleft}{\kern0pt}j{\isacharcomma}{\kern0pt}\ env{\isacharparenright}{\kern0pt}\ {\isacharequal}{\kern0pt}\ y{\isacharsemicolon}{\kern0pt}\ env\ {\isasymin}\ list{\isacharparenleft}{\kern0pt}A{\isacharparenright}{\kern0pt}{\isasymrbrakk}\isanewline
{\isasymLongrightarrow}\ x\ {\isasymin}\ y\ {\isasymlongleftrightarrow}\ A{\isacharcomma}{\kern0pt}\ env\ {\isasymTurnstile}\ {\isasymcdot}i\ {\isasymin}\ j{\isasymcdot}\isasep\isanewline\isanewline%
{\isasymlbrakk}nth{\isacharparenleft}{\kern0pt}i{\isacharcomma}{\kern0pt}\ env{\isacharparenright}{\kern0pt}\ {\isacharequal}{\kern0pt}\ x{\isacharsemicolon}{\kern0pt}\ nth{\isacharparenleft}{\kern0pt}j{\isacharcomma}{\kern0pt}\ env{\isacharparenright}{\kern0pt}\ {\isacharequal}{\kern0pt}\ y{\isacharsemicolon}{\kern0pt}\ env\ {\isasymin}\ list{\isacharparenleft}{\kern0pt}A{\isacharparenright}{\kern0pt}{\isasymrbrakk}\isanewline
{\isasymLongrightarrow}\ x\ {\isacharequal}{\kern0pt}\ y\ {\isasymlongleftrightarrow}\ A{\isacharcomma}{\kern0pt}\ env\ {\isasymTurnstile}\ {\isasymcdot}i\ {\isacharequal}{\kern0pt}\ j{\isasymcdot}\isasep\isanewline\isanewline%
env\ {\isasymin}\ list{\isacharparenleft}{\kern0pt}A{\isacharparenright}{\kern0pt}\ {\isasymLongrightarrow}\ {\isacharparenleft}{\kern0pt}A{\isacharcomma}{\kern0pt}\ env\ {\isasymTurnstile}\ {\isasymcdot}{\isasymnot}{\isacharparenleft}{\kern0pt}p\ {\isasymand}\ q{\isacharparenright}{\kern0pt}{\isasymcdot}{\isacharparenright}{\kern0pt}\ {\isasymlongleftrightarrow}\ {\isasymnot}\ {\isacharparenleft}{\kern0pt}{\isacharparenleft}{\kern0pt}A{\isacharcomma}{\kern0pt}\ env\ {\isasymTurnstile}\ p{\isacharparenright}{\kern0pt}\ {\isasymand}\isanewline%
\ \ {\isacharparenleft}{\kern0pt}A{\isacharcomma}{\kern0pt}\ env\ {\isasymTurnstile}\ q{\isacharparenright}{\kern0pt}{\isacharparenright}{\kern0pt}\isasep\isanewline\isanewline%
env\ {\isasymin}\ list{\isacharparenleft}{\kern0pt}A{\isacharparenright}{\kern0pt}\ {\isasymLongrightarrow}\ {\isacharparenleft}{\kern0pt}A{\isacharcomma}{\kern0pt}\ env\ {\isasymTurnstile}\ {\isacharparenleft}{\kern0pt}{\isasymcdot}{\isasymforall}p{\isasymcdot}{\isacharparenright}{\kern0pt}{\isacharparenright}{\kern0pt}\ {\isasymlongleftrightarrow}\ {\isacharparenleft}{\kern0pt}{\isasymforall}x{\isasymin}A{\isachardot}{\kern0pt}\ A{\isacharcomma}{\kern0pt}\ Cons{\isacharparenleft}{\kern0pt}x{\isacharcomma}{\kern0pt}\ env{\isacharparenright}{\kern0pt}\ {\isasymTurnstile}\ p{\isacharparenright}{\kern0pt}%
\end{isabelle}%
as well as the satisfaction of an arbitrary set of sentences.%
\begin{isabelle}%
A\ {\isasymTurnstile}\ {\isasymPhi}\ {\isasymequiv}\ {\isasymforall}{\isasymphi}{\isasymin}{\isasymPhi}{\isachardot}{\kern0pt}\ A{\isacharcomma}{\kern0pt}\ {\isacharbrackleft}{\kern0pt}{\isacharbrackright}{\kern0pt}\ {\isasymTurnstile}\ {\isasymphi}%
\end{isabelle}%
The internalized (viz. as elements of the set \isa{formula})
versions of the axioms are checked next against the relative statements.%
\begin{isabelle}%
Union{\isacharunderscore}{\kern0pt}ax{\isacharparenleft}{\kern0pt}{\isacharhash}{\kern0pt}{\isacharhash}{\kern0pt}A{\isacharparenright}{\kern0pt}\ {\isasymlongleftrightarrow}\ A{\isacharcomma}{\kern0pt}\ {\isacharbrackleft}{\kern0pt}{\isacharbrackright}{\kern0pt}\ {\isasymTurnstile}\ {\isasymcdot}Union\ Ax{\isasymcdot}\isasep\isanewline%
power{\isacharunderscore}{\kern0pt}ax{\isacharparenleft}{\kern0pt}{\isacharhash}{\kern0pt}{\isacharhash}{\kern0pt}A{\isacharparenright}{\kern0pt}\ {\isasymlongleftrightarrow}\ A{\isacharcomma}{\kern0pt}\ {\isacharbrackleft}{\kern0pt}{\isacharbrackright}{\kern0pt}\ {\isasymTurnstile}\ {\isasymcdot}Powerset\ Ax{\isasymcdot}\isasep\isanewline%
upair{\isacharunderscore}{\kern0pt}ax{\isacharparenleft}{\kern0pt}{\isacharhash}{\kern0pt}{\isacharhash}{\kern0pt}A{\isacharparenright}{\kern0pt}\ {\isasymlongleftrightarrow}\ A{\isacharcomma}{\kern0pt}\ {\isacharbrackleft}{\kern0pt}{\isacharbrackright}{\kern0pt}\ {\isasymTurnstile}\ {\isasymcdot}Pairing{\isasymcdot}\isasep\isanewline%
foundation{\isacharunderscore}{\kern0pt}ax{\isacharparenleft}{\kern0pt}{\isacharhash}{\kern0pt}{\isacharhash}{\kern0pt}A{\isacharparenright}{\kern0pt}\ {\isasymlongleftrightarrow}\ A{\isacharcomma}{\kern0pt}\ {\isacharbrackleft}{\kern0pt}{\isacharbrackright}{\kern0pt}\ {\isasymTurnstile}\ {\isasymcdot}Foundation{\isasymcdot}\isasep\isanewline%
extensionality{\isacharparenleft}{\kern0pt}{\isacharhash}{\kern0pt}{\isacharhash}{\kern0pt}A{\isacharparenright}{\kern0pt}\ {\isasymlongleftrightarrow}\ A{\isacharcomma}{\kern0pt}\ {\isacharbrackleft}{\kern0pt}{\isacharbrackright}{\kern0pt}\ {\isasymTurnstile}\ {\isasymcdot}Extensionality{\isasymcdot}\isasep\isanewline%
infinity{\isacharunderscore}{\kern0pt}ax{\isacharparenleft}{\kern0pt}{\isacharhash}{\kern0pt}{\isacharhash}{\kern0pt}A{\isacharparenright}{\kern0pt}\ {\isasymlongleftrightarrow}\ A{\isacharcomma}{\kern0pt}\ {\isacharbrackleft}{\kern0pt}{\isacharbrackright}{\kern0pt}\ {\isasymTurnstile}\ {\isasymcdot}Infinity{\isasymcdot}\isasep\isanewline\isanewline%
{\isasymphi}\ {\isasymin}\ formula\ {\isasymLongrightarrow}\isanewline
{\isacharparenleft}{\kern0pt}M{\isacharcomma}{\kern0pt}\ {\isacharbrackleft}{\kern0pt}{\isacharbrackright}{\kern0pt}\ {\isasymTurnstile}\ {\isasymcdot}Separation{\isacharparenleft}{\kern0pt}{\isasymphi}{\isacharparenright}{\kern0pt}{\isasymcdot}{\isacharparenright}{\kern0pt}\ {\isasymlongleftrightarrow}\isanewline
{\isacharparenleft}{\kern0pt}{\isasymforall}env{\isasymin}list{\isacharparenleft}{\kern0pt}M{\isacharparenright}{\kern0pt}{\isachardot}{\kern0pt}\isanewline
\isaindent{{\isacharparenleft}{\kern0pt}\ \ \ }arity{\isacharparenleft}{\kern0pt}{\isasymphi}{\isacharparenright}{\kern0pt}\ {\isasymle}\ {\isadigit{1}}\ {\isacharplus}{\kern0pt}\isactrlsub {\isasymomega}\ length{\isacharparenleft}{\kern0pt}env{\isacharparenright}{\kern0pt}\ {\isasymlongrightarrow}\isanewline
\ \ \ \ separation{\isacharparenleft}{\kern0pt}{\isacharhash}{\kern0pt}{\isacharhash}{\kern0pt}M{\isacharcomma}{\kern0pt}\ {\isasymlambda}x{\isachardot}{\kern0pt}\ M{\isacharcomma}{\kern0pt}\ {\isacharbrackleft}{\kern0pt}x{\isacharbrackright}{\kern0pt}\ {\isacharat}{\kern0pt}\ env\ {\isasymTurnstile}\ {\isasymphi}{\isacharparenright}{\kern0pt}{\isacharparenright}{\kern0pt}\isasep\isanewline\isanewline%
{\isasymphi}\ {\isasymin}\ formula\ {\isasymLongrightarrow}\isanewline
{\isacharparenleft}{\kern0pt}M{\isacharcomma}{\kern0pt}\ {\isacharbrackleft}{\kern0pt}{\isacharbrackright}{\kern0pt}\ {\isasymTurnstile}\ {\isasymcdot}Replacement{\isacharparenleft}{\kern0pt}{\isasymphi}{\isacharparenright}{\kern0pt}{\isasymcdot}{\isacharparenright}{\kern0pt}\ {\isasymlongleftrightarrow}\ {\isacharparenleft}{\kern0pt}{\isasymforall}env{\isachardot}{\kern0pt}\ replacement{\isacharunderscore}{\kern0pt}assm{\isacharparenleft}{\kern0pt}M{\isacharcomma}{\kern0pt}\ env{\isacharcomma}{\kern0pt}\ {\isasymphi}{\isacharparenright}{\kern0pt}{\isacharparenright}\isanewline\isanewline%
choice{\isacharunderscore}{\kern0pt}ax{\isacharparenleft}{\kern0pt}{\isacharhash}{\kern0pt}{\isacharhash}{\kern0pt}A{\isacharparenright}{\kern0pt}\ {\isasymlongleftrightarrow}\ A{\isacharcomma}{\kern0pt}\ {\isacharbrackleft}{\kern0pt}{\isacharbrackright}{\kern0pt}\ {\isasymTurnstile}\ {\isasymcdot}AC{\isasymcdot}%
\end{isabelle}%

Finally, the axiom sets are defined as follows.

\begin{isabelle}%
ZF{\isacharunderscore}{\kern0pt}fin\ {\isasymequiv}\isanewline
{\isacharbraceleft}{\kern0pt}{\isasymcdot}Extensionality{\isasymcdot}{\isacharcomma}{\kern0pt}\ {\isasymcdot}Foundation{\isasymcdot}{\isacharcomma}{\kern0pt}\ {\isasymcdot}Pairing{\isasymcdot}{\isacharcomma}{\kern0pt}\ {\isasymcdot}Union\ Ax{\isasymcdot}{\isacharcomma}{\kern0pt}\ {\isasymcdot}Infinity{\isasymcdot}{\isacharcomma}{\kern0pt}\isanewline
\isaindent{{\isacharbraceleft}{\kern0pt}}{\isasymcdot}Powerset\ Ax{\isasymcdot}{\isacharbraceright}{\kern0pt}\isasep\isanewline\isanewline%
ZF{\isacharunderscore}{\kern0pt}schemes\ {\isasymequiv}\isanewline
{\isacharbraceleft}{\kern0pt}{\isasymcdot}Separation{\isacharparenleft}{\kern0pt}p{\isacharparenright}{\kern0pt}{\isasymcdot}\ {\isachardot}{\kern0pt}\ p\ {\isasymin}\ formula{\isacharbraceright}{\kern0pt}\ {\isasymunion}\ {\isacharbraceleft}{\kern0pt}{\isasymcdot}Replacement{\isacharparenleft}{\kern0pt}p{\isacharparenright}{\kern0pt}{\isasymcdot}\ {\isachardot}{\kern0pt}\ p\ {\isasymin}\ formula{\isacharbraceright}{\kern0pt}\isasep\isanewline\isanewline%
{\isasymcdot}Z{\isasymcdot}\ {\isasymequiv}\ ZF{\isacharunderscore}{\kern0pt}fin\ {\isasymunion}\ {\isacharbraceleft}{\kern0pt}{\isasymcdot}Separation{\isacharparenleft}{\kern0pt}p{\isacharparenright}{\kern0pt}{\isasymcdot}\ {\isachardot}{\kern0pt}\ p\ {\isasymin}\ formula{\isacharbraceright}{\kern0pt}\isasep\isanewline%
ZC\ {\isasymequiv}\ {\isasymcdot}Z{\isasymcdot}\ {\isasymunion}\ {\isacharbraceleft}{\kern0pt}{\isasymcdot}AC{\isasymcdot}{\isacharbraceright}{\kern0pt}\isasep\isanewline%
ZF\ {\isasymequiv}\ ZF{\isacharunderscore}{\kern0pt}schemes\ {\isasymunion}\ ZF{\isacharunderscore}{\kern0pt}fin\isasep\isanewline%
ZFC\ {\isasymequiv}\ ZF\ {\isasymunion}\ {\isacharbraceleft}{\kern0pt}{\isasymcdot}AC{\isasymcdot}{\isacharbraceright}{\kern0pt}%
\end{isabelle}%

\subsection{Relativization of infinitary arithmetic\label{sec:relative-arith}%
}
In order to state the defining property of the relative
equipotency relation, we work under the assumptions of the
locale \isa{M{\isacharunderscore}{\kern0pt}cardinals}. They comprise a finite set
of instances of Separation and Replacement to prove
closure properties of the transitive class \isa{M}.%
\begin{isabelle}
\isacommand{lemma}\isamarkupfalse%
\ {\isacharparenleft}{\kern0pt}\isakeyword{in}\ M{\isacharunderscore}{\kern0pt}cardinals{\isacharparenright}{\kern0pt}\ eqpoll{\isacharunderscore}{\kern0pt}def{\isacharprime}{\kern0pt}{\isacharcolon}{\kern0pt}\isanewline
\ \ \isakeyword{assumes}\ {\isachardoublequoteopen}M{\isacharparenleft}{\kern0pt}A{\isacharparenright}{\kern0pt}{\isachardoublequoteclose}\ {\isachardoublequoteopen}M{\isacharparenleft}{\kern0pt}B{\isacharparenright}{\kern0pt}{\isachardoublequoteclose}\ \isakeyword{shows}\ {\isachardoublequoteopen}A\ {\isasymapprox}\isactrlbsup M\isactrlesup \ B\ {\isasymlongleftrightarrow}\ {\isacharparenleft}{\kern0pt}{\isasymexists}f{\isacharbrackleft}{\kern0pt}M{\isacharbrackright}{\kern0pt}{\isachardot}{\kern0pt}\ f\ {\isasymin}\ bij{\isacharparenleft}{\kern0pt}A{\isacharcomma}{\kern0pt}B{\isacharparenright}{\kern0pt}{\isacharparenright}{\kern0pt}{\isachardoublequoteclose}
\end{isabelle}

Below, $\mu$ denotes the minimum operator on the ordinals.%
\begin{isabelle}
  \isacommand{lemma}\isamarkupfalse%
\ cardinalities{\isacharunderscore}{\kern0pt}defs{\isacharcolon}{\kern0pt}\isanewline
\ \ \isakeyword{fixes}\ M{\isacharcolon}{\kern0pt}{\isacharcolon}{\kern0pt}{\isachardoublequoteopen}i{\isasymRightarrow}o{\isachardoublequoteclose}\isanewline
\ \ \isakeyword{shows}\isanewline
\ \ \ \ {\isachardoublequoteopen}{\isacharbar}{\kern0pt}A{\isacharbar}{\kern0pt}\isactrlbsup M\isactrlesup \ {\isasymequiv}\ {\isasymmu}\ i{\isachardot}{\kern0pt}\ M{\isacharparenleft}{\kern0pt}i{\isacharparenright}{\kern0pt}\ {\isasymand}\ i\ {\isasymapprox}\isactrlbsup M\isactrlesup \ A{\isachardoublequoteclose}\isanewline
\ \ \ \ {\isachardoublequoteopen}Card\isactrlbsup M\isactrlesup {\isacharparenleft}{\kern0pt}{\isasymalpha}{\isacharparenright}{\kern0pt}\ {\isasymequiv}\ {\isasymalpha}\ {\isacharequal}{\kern0pt}\ {\isacharbar}{\kern0pt}{\isasymalpha}{\isacharbar}{\kern0pt}\isactrlbsup M\isactrlesup {\isachardoublequoteclose}\isanewline
\ \ \ \ {\isachardoublequoteopen}{\isasymkappa}\isactrlbsup {\isasymup}{\isasymnu}{\isacharcomma}{\kern0pt}M\isactrlesup \ {\isasymequiv}\ {\isacharbar}{\kern0pt}{\isasymnu}\ {\isasymrightarrow}\isactrlbsup M\isactrlesup \ {\isasymkappa}{\isacharbar}{\kern0pt}\isactrlbsup M\isactrlesup {\isachardoublequoteclose}\isanewline
\ \ \ \ {\isachardoublequoteopen}{\isacharparenleft}{\kern0pt}{\isasymkappa}\isactrlsup {\isacharplus}{\kern0pt}{\isacharparenright}{\kern0pt}\isactrlbsup M\isactrlesup \ {\isasymequiv}\ {\isasymmu}\ x{\isachardot}{\kern0pt}\ M{\isacharparenleft}{\kern0pt}x{\isacharparenright}{\kern0pt}\ {\isasymand}\ Card\isactrlbsup M\isactrlesup {\isacharparenleft}{\kern0pt}x{\isacharparenright}{\kern0pt}\ {\isasymand}\ {\isasymkappa}\ {\isacharless}{\kern0pt}\ x{\isachardoublequoteclose}
\end{isabelle}
Analogous to the previous Lemma
\isa{eqpoll{\isacharunderscore}{\kern0pt}def{\isacharprime}{\kern0pt}},
the next lemma holds under
the assumptions of the locale \isa{M{\isacharunderscore}{\kern0pt}aleph}. The axiom instances
included are sufficient to state and prove the defining
properties of the relativized \isa{Aleph} function
(in particular, the required ability to perform transfinite recursions).%
\begin{isabelle}%
\isacommand{context}\isamarkupfalse%
\ M{\isacharunderscore}{\kern0pt}aleph\isanewline
\isakeyword{begin}%
\isanewline
\isanewline
{\isasymaleph}\isactrlbsub {\isadigit{0}}\isactrlesub \isactrlbsup M\isactrlesup \ {\isacharequal}{\kern0pt}\ {\isasymomega}\isasep\isanewline%
{\isasymlbrakk}Ord{\isacharparenleft}{\kern0pt}{\isasymalpha}{\isacharparenright}{\kern0pt}{\isacharsemicolon}{\kern0pt}\ M{\isacharparenleft}{\kern0pt}{\isasymalpha}{\isacharparenright}{\kern0pt}{\isasymrbrakk}\ {\isasymLongrightarrow}\ {\isasymaleph}\isactrlbsub succ{\isacharparenleft}{\kern0pt}{\isasymalpha}{\isacharparenright}{\kern0pt}\isactrlesub \isactrlbsup M\isactrlesup \ {\isacharequal}{\kern0pt}\ {\isacharparenleft}{\kern0pt}{\isasymaleph}\isactrlbsub {\isasymalpha}\isactrlesub \isactrlbsup M\isactrlesup \isactrlsup {\isacharplus}{\kern0pt}{\isacharparenright}{\kern0pt}\isactrlbsup M\isactrlesup \isasep\isanewline%
{\isasymlbrakk}Limit{\isacharparenleft}{\kern0pt}{\isasymalpha}{\isacharparenright}{\kern0pt}{\isacharsemicolon}{\kern0pt}\ M{\isacharparenleft}{\kern0pt}{\isasymalpha}{\isacharparenright}{\kern0pt}{\isasymrbrakk}\ {\isasymLongrightarrow}\ {\isasymaleph}\isactrlbsub {\isasymalpha}\isactrlesub \isactrlbsup M\isactrlesup \ {\isacharequal}{\kern0pt}\ {\isacharparenleft}{\kern0pt}{\isasymUnion}j{\isasymin}{\isasymalpha}{\isachardot}{\kern0pt}\ {\isasymaleph}\isactrlbsub j\isactrlesub \isactrlbsup M\isactrlesup {\isacharparenright}{\kern0pt}%
\end{isabelle}%
\isacommand{end}\isamarkupfalse%
\ %
\isamarkupcmt{\isa{M{\isacharunderscore}{\kern0pt}aleph}%
}
\begin{isabelle}
\isacommand{lemma}\isamarkupfalse%
\ ContHyp{\isacharunderscore}{\kern0pt}rel{\isacharunderscore}{\kern0pt}def{\isacharprime}{\kern0pt}{\isacharcolon}{\kern0pt}\isanewline
\ \ \isakeyword{fixes}\ N{\isacharcolon}{\kern0pt}{\isacharcolon}{\kern0pt}{\isachardoublequoteopen}i{\isasymRightarrow}o{\isachardoublequoteclose}\isanewline
\ \ \isakeyword{shows}\isanewline
\ \ \ \ {\isachardoublequoteopen}CH\isactrlbsup N\isactrlesup \ {\isasymequiv}\ {\isasymaleph}\isactrlbsub {\isadigit{1}}\isactrlesub \isactrlbsup N\isactrlesup \ {\isacharequal}{\kern0pt}\ {\isadigit{2}}\isactrlbsup {\isasymup}{\isasymaleph}\isactrlbsub {\isadigit{0}}\isactrlesub \isactrlbsup N\isactrlesup {\isacharcomma}{\kern0pt}N\isactrlesup {\isachardoublequoteclose}
\end{isabelle}

Under appropriate hypotheses (this time, from the locale \isa{M{\isacharunderscore}{\kern0pt}ZF{\isacharunderscore}{\kern0pt}library}),
   \isa{CH\isactrlbsup M\isactrlesup } is equivalent to its fully relational version \isa{is{\isacharunderscore}{\kern0pt}ContHyp}.
    As a sanity check, we see that if the transitive class is indeed \isa{{\isasymV}},
    we recover the original $\CH$.%
\begin{isabelle}%
M{\isacharunderscore}{\kern0pt}ZF{\isacharunderscore}{\kern0pt}library{\isacharparenleft}{\kern0pt}M{\isacharparenright}{\kern0pt}\ {\isasymLongrightarrow}\ is{\isacharunderscore}{\kern0pt}ContHyp{\isacharparenleft}{\kern0pt}M{\isacharparenright}{\kern0pt}\ {\isasymlongleftrightarrow}\ CH\isactrlbsup M\isactrlesup \isasep\isanewline%
is{\isacharunderscore}{\kern0pt}ContHyp{\isacharparenleft}{\kern0pt}{\isasymV}{\isacharparenright}{\kern0pt}\ {\isasymlongleftrightarrow}\ {\isasymaleph}\isactrlbsub {\isadigit{1}}\isactrlesub \ {\isacharequal}{\kern0pt}\ {\isadigit{2}}\isactrlbsup {\isasymup}{\isasymaleph}\isactrlbsub {\isadigit{0}}\isactrlesub \isactrlesup %
\end{isabelle}%
In turn, the fully relational version evaluated on a nonempty
transitive \isa{A} is equivalent to the satisfaction of the
first-order formula \isa{{\isasymcdot}CH{\isasymcdot}} (since it
actually is a sentence, it does not depend on \isa{env}, which
appears only because the definition of $\models$ requires that argument).%
\begin{isabelle}%
{\isasymlbrakk}env\ {\isasymin}\ list{\isacharparenleft}{\kern0pt}A{\isacharparenright}{\kern0pt}{\isacharsemicolon}{\kern0pt}\ {\isadigit{0}}\ {\isasymin}\ A{\isasymrbrakk}\ {\isasymLongrightarrow}\ is{\isacharunderscore}{\kern0pt}ContHyp{\isacharparenleft}{\kern0pt}{\isacharhash}{\kern0pt}{\isacharhash}{\kern0pt}A{\isacharparenright}{\kern0pt}\ {\isasymlongleftrightarrow}\ A{\isacharcomma}{\kern0pt}\ env\ {\isasymTurnstile}\ {\isasymcdot}CH{\isasymcdot}%
\end{isabelle}%

\section{Discipline for relativization}
\label{sec:discipline-relativization}

As we said in Sec.~\ref{sec:relat-vers-non-absol}, in order to force
$\CH$ and its negation we depended on having relativized versions of
cardinals, Alephs, etc. It was clear for us that our efforts would be
more efficient if we set up a discipline for relativizing sets (terms
of type $\tyi$) and predicates/relations (terms of type $\tyo$).

Paulson only had, for each set, the relational version. It seemed
clearer to us to have a functional version of the relativized concept.
Going back to our example in \ref{sec:tools-relativization}, for the
concept $\isa{cardinal}::\tyi \fun \tyi$ we want its relative
version
$\isa{cardinal{\uscore}rel}::(\tyi \fun \tyo) \fun \tyi \fun\tyi$
and the relational version of the latter
$\isa{is{\uscore}cardinal}::(\tyi \fun \tyo) \fun \tyi \fun \tyi
\fun \tyo$.

Our first attempt of defining a discipline was inspired by
mathematical considerations: if we might prove that
$\isa{is{\uscore}cardinal}$ is functional and also prove the
existence of a witness $\isa{c}$ such that $\isa{M(c)}$ and
$\isa{is{\uscore}cardinal(M,x,c)}$ then
$\isa{cardinal{\uscore}rel(M,x)}$ can be obtained by the operator of
definite descriptions.

Soon we realized that resorting to definite descriptions was needed
only for the most primitive concepts. In fact, once we have a
relativized concept, we can use it to define other relativizations.
For instance, $\isa{cardinal{\uscore}rel}$ depends on having
relative versions of $\isa{bij}$. Instead of relationalizing
$\isa{bij}$ to get $\isa{is{\uscore}bij}$ and then prove
uniqueness and existence of a witness, we define
$\isa{bij{\uscore}rel}$ using $\isa{inj{\uscore}rel}$ and
$\isa{surj{\uscore}rel}$.

\section{Recursions in cofinality and the Delta System Lemma}\label{sec:recursions-cofinality}

As we mentioned near the end of
Section~\ref{sec:aims-formalization-planning}, we decided to minimize
the requirements being formalized in order to achieve our immediate
goal. In particular, the treatment of cofinality in the companion
project \cite{Delta_System_Lemma-AFP} was left behind.

We already observed that well-founded, and in particular transfinite,
recursion is not easily dealt with in Isabelle/ZF. Nevertheless, and
mainly as a curiosity, we found out that only one recursive
construction is needed for the development of the basic theory of
cofinality (as in \cite[Sect.~I.13]{kunen2011set}), which is used in
the proof of the following “factorization” lemma:

\begin{lemma}
  Let $\del,\ga\in\Ord$ and assume $f:\del\to\ga$ is cofinal.  There exists
  a strictly increasing $g:\cf(\ga)\to \del$ such that $f\circ g$ is
  strictly increasing and cofinal in $\ga$. Moreover, if $f$ is
  strictly increasing, then $g$ must also be cofinal.
\end{lemma}

It turns out that the rest of the basic results on cofinality (namely,
idempotence of $\cf$, that regular ordinals are cardinals, the
cofinality of Alephs, König's Theorem) follow easily from the previous
Lemma by “algebraic” reasoning only.
We expect the relativization of these
results to be straightforward.

The AFP entry \cite{Delta_System_Lemma-AFP} also includes the
formalization of the (absolute) Delta System Lemma (DSL). Formalizing its
proof was rather straightforward, once the many prerequisites were
taken care of. Some of those were really basic, for instance:
\begin{enumerate}
\item \label{item:1}$\omega$ injects into every infinite set;
\item \label{item:2}surjective images of countable sets are countable;
\item \label{item:3}the union over a countable index set $J$ of a family $X :: \tyi
  \fun \tyi$ of countable sets is countable.
\end{enumerate}
It was also convenient to isolate the relevant recursive construction
principle (\isatt{bounded{\uscore}cardinal{\uscore}selection}) that
appears in the proof of DSL, which was also useful for showing
Item~\ref{item:1}:
\begin{lemma}\label{lem:bdd-card-selection}
  Assume $F$ is nonempty,
  $\gamma$ is a cardinal, and $Q$ is a binary predicate over $F$ satisfying 
  \[
    \forall Y \sbq F.\ |Y| < \gamma \implies \exists a\in F.\ \forall
    y\in Y.\ Q(y,a).
  \]
  Then there exists  $S:\gamma \to F$ such that
  $Q(S(\alpha),S(\beta))$ for all $\alpha<\beta<\gamma$.
\end{lemma}

Concerning the relativization of the proof of DSL for its use in
forcing, it required inserting proofs that all the relevant objects
lied in the model $M$; this was only tedious.
A bit more effort was required at the point where Item~\ref{item:3}
was used, because it involved an application of Replacement (being $X$
a class function); relativizing
Lemma~\ref{lem:bdd-card-selection} also required some work, because of the
recursion.

\section{Axioms of Isabelle/ZF}
\label{appendix:axioms}

In this appendix we list the complete set of axioms of Isabelle's
metatheory and logic.

\subsection{The metatheory $\Meta$}
Below, the \isatt{PROP} operator corresponds to the
injection $[\cdot]$ of Section~\ref{sec:isabelle-metalogic-meta}.
\begin{isabelle}
Pure.abstract\_rule: (⋀x. ?f(x) ≡ ?g(x)) ⟹ λx. ?f(x) ≡ λx. ?g(x)\isanewline
Pure.combination: ?f ≡ ?g ⟹ ?x ≡ ?y ⟹ ?f(?x) ≡ ?g(?y)\isanewline
Pure.equal\_elim: PROP ?A ≡ PROP ?B ⟹ PROP ?A ⟹ PROP ?B\isanewline
Pure.equal\_intr: (PROP ?A ⟹ PROP ?B) ⟹ (PROP ?B ⟹ PROP ?A) ⟹ \isanewline
\ \ \ \ \ \ \ \ \ \ PROP ?A ≡ PROP ?B\isanewline
Pure.reflexive: ?x ≡ ?x\isanewline
Pure.symmetric: ?x ≡ ?y ⟹ ?y ≡ ?x\isanewline
Pure.transitive: ?x ≡ ?y ⟹ ?y ≡ ?z ⟹ ?x ≡ ?z
\end{isabelle}

\subsection{IFOL and FOL}
In the axioms \isa{refl}, \isa{subst}, \isa{allI}, \isa{spec},
\isa{exE}, \isa{exI}, and \isa{eq\_reflection} there is a constraint
on the \emph{type} of the variables \isa{a, b, x}: It should be in the class \isa{term\_class}.
\begin{isabelle}
IFOL.FalseE: ⋀P. False ⟹ P\isanewline
IFOL.refl:  (⋀a. a = a)\isanewline
IFOL.subst: (⋀a b P. a = b ⟹ P(a) ⟹ P(b))\isanewline
IFOL.allI:  (⋀P. (⋀x. P(x)) ⟹ ∀x. P(x))\isanewline
IFOL.spec:  (⋀P x. ∀x. P(x) ⟹ P(x))\isanewline
IFOL.exE:   (⋀P R. ∃x. P(x) ⟹ (⋀x. P(x) ⟹ R) ⟹ R)\isanewline
IFOL.exI:   (⋀P x. P(x) ⟹ ∃x. P(x))\isanewline
IFOL.conjI: ⋀P Q. P ⟹ Q ⟹ P ∧ Q\isanewline
IFOL.conjunct1: ⋀P Q. P ∧ Q ⟹ P\isanewline
IFOL.conjunct2: ⋀P Q. P ∧ Q ⟹ Q\isanewline
IFOL.disjE: ⋀P Q R. P ∨ Q ⟹ (P ⟹ R) ⟹ (Q ⟹ R) ⟹ R\isanewline
IFOL.disjI1: ⋀P Q. P ⟹ P ∨ Q\isanewline
IFOL.disjI2: ⋀P Q. Q ⟹ P ∨ Q\isanewline
IFOL.eq\_reflection: (⋀x y. x = y ⟹ x ≡ y)\isanewline
IFOL.iff\_reflection: ⋀P Q. P ⟷ Q ⟹ P ≡ Q\isanewline
IFOL.impI: ⋀P Q. (P ⟹ Q) ⟹ P ⟶ Q\isanewline
IFOL.mp: ⋀P Q. P ⟶ Q ⟹ P ⟹ Q\isanewline%
FOL.classical: ⋀P. (¬ P ⟹ P) ⟹ P
\end{isabelle}

\subsection{ZF\_Base}
The following symbols are introduced in this theory:
\begin{isabelle}
axiomatization\isanewline
\ \ \ \      mem :: "[i, i] ⇒ o"  (infixl ‹∈› 50)  \isamarkupcmt{membership relation}\isanewline
  and zero :: "i"  (‹0›)  \isamarkupcmt{the empty set}\isanewline
  and Pow :: "i ⇒ i"  \isamarkupcmt{power sets}\isanewline
  and Inf :: "i"  \isamarkupcmt{infinite set}\isanewline
  and Union :: "i ⇒ i"  (‹⋃\_› [90] 90)\isanewline
  and PrimReplace :: "[i, [i, i] ⇒ o] ⇒ i"
\end{isabelle}
\noindent After the definitions of $\notin$, $\subseteq$, $\isa{succ}$,
and relative quantifications are presented, the following axioms are postulated:
\begin{isabelle}
ZF\_Base.Pow\_iff: ⋀A B. A ∈ Pow(B) ⟷ A ⊆ B\isanewline
ZF\_Base.Union\_iff: ⋀A C. A ∈ ⋃C ⟷ (∃B∈C. A ∈ B)\isanewline
ZF\_Base.extension: ⋀A B. A = B ⟷ A ⊆ B ∧ B ⊆ A\isanewline
ZF\_Base.foundation: ⋀A. A = 0 ∨ (∃x∈A. ∀y∈x. y ∉ A)\isanewline
ZF\_Base.infinity: 0 ∈ Inf ∧ (∀y∈Inf. succ(y) ∈ Inf)\isanewline
ZF\_Base.replacement: ⋀A P b. ∀x∈A. ∀y z. P(x, y) ∧ P(x, z) ⟶ y = z \isanewline
\ \ \ \ \ \ \ \ \ ⟹ b ∈ PrimReplace(A, P) ⟷ (∃x∈A. P(x, b))
\end{isabelle}

\subsection{AC}

The theory \theory{AC} is only imported in the theory
\theory{Absolute\_Versions}; none of the main results
depends on $\AC$. The latter theory
shows that some absolute results can be obtained from the
relativized versions on $\mathcal{V}$.

\begin{isabelle}
AC.AC: ⋀a A B. a ∈ A ⟹ (⋀x. x ∈ A ⟹ ∃y. y ∈ B(x)) ⟹\isanewline
  \ \ \ \ \ \ ∃z. z ∈ Pi(A, B)
\end{isabelle}

\section{Lambda replacements}\label{sec:lambda-replacements}

The development of the locale structure of the project was a dynamical
process. As further properties of closure of the ground $M$ were
required, we gathered the relevant instances of Separation and
Replacement into a new locale (always assuming a class model, for
added generality), and proceeded to apply them to those closure proofs.

This procedure lead to a steady grow in the number of interpretation
obligations and therefore, of formula synthesis (since the two axiom
schemes were postulated using codes for formulas). That number would
easily surpass the hundred, and the automatic tools at our disposal
for that task were rudimentary (as discussed in
Section~\ref{sec:bureaucracy-scale-factors}).

Facing this situation, we decided that we needed some sort of
\emph{compositionality} in order to obtain new instances from the ones
already proved: Having Replacement for class functions $F$ and $G$
does not entail immediately replacement under $F\circ G$ (unless you
use one further instance of Separation, and the net gain is zero). The
solution was to postulate for the relevant $F$s, instead of
replacement through $x\mapsto F(x)$, a \emph{lambda replacement}
through $x\mapsto \lb x,F(x)\rb$. The name “lambda” corresponds to the
fact that this type of replacement is equivalent to closure under
$(\lambda x\in A.\ F(x)) \defi \{ \lb x,F(x)\rb : x\in A \}$ for every
$A\in M$.

Now, a fixed set of six replacements and one separation (apart from
those in \locale{M{\uscore}basic}, which also assumes the Powerset
Axiom for the class $M$) is sufficient to obtain the lambda
replacement under $x\mapsto \lb x,F(G(x))\rb$ given those for $F$ and
$G$. To obtain compositions with binary class functions $H$, it is
enough to assume the lambda replacement
$x\mapsto \lb x,H(\mathit{fst}(x),\mathit{snd}(x)))\rb$. We summarize
the assumptions in Table~\ref{tab:m-repl-instances}.

\newcommand{\lamRepl}[2][x]{#1 \mapsto\langle #1,#2\rangle}
\begin{table}[!h]
\centering
\begin{threeparttable}
\begin{tabular}{r<{\stepcounter{LamReplCount}\theLamReplCount.} >{\hspace{1ex}}l @{\hspace{0.8em}} p{16.3em}}
  \toprule
  \multicolumn{1}{r}{No.} & Name & Instance  \\
  \midrule
  \multicolumn{3}{@{}l}{\hspace{0.3em}{\textit{{Replacement Instances}}}}\\
  & \isa{lam{\uscore}replacement{\uscore}fst} & $\lamRepl{\mathit{fst}(x)}$ \\
  & \isa{lam{\uscore}replacement{\uscore}snd} & $\lamRepl{\mathit{snd}(x)}$ \\
  & \isa{lam{\uscore}replacement{\uscore}Union} & $\lamRepl{\bigcup(x)}$ \\
  & \isa{lam{\uscore}replacement{\uscore}Image} & $\lamRepl{\mathit{fst}(x)``\mathit{snd}(x)}$ \\
  & \isa{lam{\uscore}replacement{\uscore}middle{\uscore}del} &
        $\lamRepl{\langle \mathit{fst}(\mathit{fst}(x)),\mathit{snd}(\mathit{snd}(x)) \rangle}$ \\
  & \isa{lam{\uscore}replacement{\uscore}prodRepl} &
  $x \mapsto \bigl\langle \mathit{snd}(\mathit{fst}(x)),\bigr.$
  \newline
  $\bigl.\phantom{x \mapsto \langle \mathit{snd}}\langle \mathit{fst}(\mathit{fst}(x)),\mathit{snd}(\mathit{snd}(x))\rangle \bigr\rangle$\\
  \midrule
  \multicolumn{3}{@{}l}{\hspace{0.3em}{\textit{{Separation Instances}}}}\\
  & \isa{middle{\uscore}separation} & $\mathit{snd}(\mathit{fst}(x))=\mathit{fst}(\mathit{snd}(x))$ \\
  & \isa{separation{\uscore}fst{\uscore}in{\uscore}snd} & $\mathit{fst}(\mathit{snd}(x)) \in \mathit{snd}(\mathit{snd}(x))$\\
  \bottomrule
\end{tabular}
\caption{Replacement and Separation instances of the locale \isa{M{\uscore}replacement}.}
\label{tab:m-repl-instances}
\end{threeparttable}
\end{table}

\section{22 replacement instances to rule them all}
\label{sec:repl-instances-appendix}

In Table~\ref{tab:instances1} we show the name of the fourteen
formulas involved in the twenty two instances of replacement needed in
our mechanization. The formulas marked with ($\dagger$) are needed
twice: one by themselves and the other as the argument for
$\isa{ground{\uscore}repl{\uscore}fm}$. The
$\isa{ground{\uscore}repl{\uscore}fm}$ function maps $\phi$ to $\psi$
as described in Sect.~\ref{sec:repl-instances}. These eight instances form
the set \isa{instances3{\uscore}fms}.

\newcommand{\groundRepl}{\ensuremath{{}^\dagger}}
\newcommand{\replInstSet}[1]{\multicolumn{3}{@{}l}{\hspace{0.3em}{$\mathbf{\mathit{#1}}$}}}

\begin{table}[!h]
\centering
\begin{threeparttable}
\begin{tabular}{r<{\stepcounter{replInstCount}\thereplInstCount.} >{\hspace{2pt}}l @{\hspace{1ex}} p{16.8em}}
  \toprule
  \multicolumn{1}{r}{No.} & Formula name & Comment \\
  \midrule
  \replInstSet{instances1{\uscore}fms}\\
  & \isa{eclose{\uscore}closed{\uscore}fm} \groundRepl & Closure under iteration of $X\mapsto\union X$. \\
  & \isa{eclose{\uscore}abs{\uscore}fm} \groundRepl & Absoluteness of the previous construction.\\
  & \isa{wfrec{\uscore}rank{\uscore}fm} \groundRepl &  For $\in$-rank.\\
  & \isa{transrec{\uscore}VFrom{\uscore}fm} \groundRepl & For rank initial segments.\\
  \midrule
  \replInstSet{instances2{\uscore}fms}\\
  & \isa{wfrec{\uscore}ordertype{\uscore}fm} \groundRepl & Well-founded recursion for the construction of ordertypes. \\
  & \isa{omap{\uscore}replacement{\uscore}fm} \groundRepl & Auxiliary instance for the definition of ordertypes. \\
  & \isa{ordtype{\uscore}replacement{\uscore}fm} \groundRepl& Replacement through the ordertype function, for Hartogs' Theorem.\\
  & \isa{wfrec{\uscore}Aleph{\uscore}fm} \groundRepl& Well-founded recursion to define Aleph.\\
  \midrule
  \replInstSet{instances{\uscore}ground{\uscore}fms}\\
  & \isa{wfrec{\uscore}Hcheck{\uscore}fm} & Well-founded recursion to define check.\\
  & \isa{wfrec{\uscore}Hfrc{\uscore}at{\uscore}fm}. & Well-founded recursion for the definition of forcing for atomic formulas.\\
  & \isa{lam{\uscore}replacement{\uscore}check{\uscore}fm} & Replacement through $x\mapsto \lb x,\check{x}\rb$, for $\punto{G}$.\\
  \midrule
  \replInstSet{instances{\uscore}ground{\uscore}notCH{\uscore}fms}\\
  &   \isa{rec{\uscore}constr{\uscore}fm} &
  Recursive construction of sets using a choice function.\\
  & \isa{rec{\uscore}constr{\uscore}abs{\uscore}fm} &
  Absoluteness of the previous construction.\\
  \midrule
  \replInstSet{instances{\uscore}ground{\uscore}CH{\uscore}fms}\\
  & \isa{dc{\uscore}abs{\uscore}fm} &  Absoluteness of the recursive construction in the proof of
  Dependent Choice from $\AC$. \\
  \midrule
  \replInstSet{instances3{\uscore}fms}\\
  \multicolumn{1}{r}{15-22.} & $\isa{ground{\uscore}repl{\uscore}fm}(\phi)$ & one for each formula $\phi$ marked with $\dagger$ \\
  \bottomrule
\end{tabular}
\caption{Replacement Instances used in our mechanization.}
\label{tab:instances1}
\end{threeparttable}
\end{table}

\end{document}